\documentclass[a4paper]{amsart} 
\usepackage[ascii]{inputenc} 

\usepackage{amssymb}
\usepackage{mathrsfs}
\usepackage[hidelinks]{hyperref}

\usepackage{xcolor}

\usepackage[shortlabels]{enumitem}
\setlist[enumerate,1]{label={(\Alph*)}}
\setlist[enumerate,2]{label={(\alph*)}}
\setlist[enumerate,3]{label={$\bullet_{\arabic*}$}}

\newenvironment{PROOF}[2][\proofname.]
{\begin{proof}[#1]\renewcommand{\qedsymbol}{\textsquare$_{\rm #2}$}}
{\end{proof}}

\newtheorem{theorem}{Theorem}[section]

\newtheorem{claim}[theorem]{Claim}

\newtheorem{observation}[theorem]{Observation}

\theoremstyle{definition}

\newtheorem{definition}[theorem]{Definition}
\newtheorem{discussion}[theorem]{Discussion}
\newtheorem{example}[theorem]{Example}

\newtheorem{fact}[theorem]{Fact}
\newtheorem{hypothesis}[theorem]{Hypothesis}

\theoremstyle{remark}

\newtheorem{notation}[theorem]{Notation}

\newtheorem{remark}[theorem]{Remark}

\newcommand{\Card}{\mathrm{Card}}

\newcommand{\plus}{\mathsf{plus}}
\newcommand{\paut}{\mathrm{paut}}

\newcommand{\bfTR}{\mathbf{TR}}
\newcommand{\TR}{\mathrm{TR}}

\newcommand{\AC}{\mathsf{AC}}

\newcommand{\DC}{\mathsf{DC}}

\newcommand{\ZF}{\mathsf{ZF}}
\newcommand{\ZFC}{\mathsf{ZFC}}

\newcommand{\Cohen}{\mathsf{Cohen}}

\newcommand{\aut}{\mathrm{aut}}
\newcommand{\cf}{\mathrm{cf}}

\newcommand{\dom}{\mathrm{dom}}
\newcommand{\HOD}{\mathrm{HOD}}

\newcommand{\id}{\mathrm{id}}
\newcommand{\inc}{\mathrm{inc}}

\newcommand{\Ord}{\mathrm{Ord}}
\newcommand{\otp}{\mathrm{otp}}
\newcommand{\pre}{\mathrm{pre}}

\newcommand{\rang}{\mathrm{rang}}

\newcommand{\at}{\mathrm{at}}

\newcommand{\pr}{\mathrm{pr}}

\newcommand{\supp}{\mathrm{supp}}

\newcommand{\tr}{\mathrm{tr}}

\newcommand{\val}{\mathrm{val}}

\newcommand{\Reg}{\mathrm{Reg}}

\newcommand{\bfB}{\mathbf{B}}

\newcommand{\bfG}{\mathbf{G}}

\newcommand{\bfK}{\mathbf{K}}
\newcommand{\bfL}{\mathbf{L}}
\newcommand{\bfM}{\mathbf{M}}
\newcommand{\bfN}{\mathbf{N}}
\newcommand{\bfO}{\mathbf{O}}
\newcommand{\bfP}{\mathbf{P}}

\newcommand{\bfQ}{\mathbf{Q}}

\newcommand{\bfV}{\mathbf{V}}

\newcommand{\bfY}{\mathbf{Y}}

\newcommand{\bfd}{\mathbf{d}}

\newcommand{\bfh}{\mathbf{h}}
\newcommand{\bfi}{\mathbf{i}}

\newcommand{\bfm}{\mathbf{m}}
\newcommand{\bfn}{\mathbf{n}}
\newcommand{\bfo}{\mathbf{o}}
\newcommand{\bfp}{\mathbf{p}}
\newcommand{\bfq}{\mathbf{q}}
\newcommand{\bfr}{\mathbf{r}}

\newcommand{\bfx}{\mathbf{x}}
\newcommand{\bfy}{\mathbf{y}}

\newcommand{\bbL}{\mathbb{L}}

\newcommand{\bbP}{\mathbb{P}}
\newcommand{\bbQ}{\mathbb{Q}}
\newcommand{\bbR}{\mathbb{R}}

\newcommand{\mn}{\medskip\noindent}
\newcommand{\sn}{\smallskip\noindent}
\newcommand{\bn}{\bigskip\noindent}

\newcommand{\cP}{\mathscr{P}}

\newcommand{\clA}{\mathcal{A}}
\newcommand{\clB}{\mathcal{B}}

\newcommand{\clE}{\mathcal{E}}
\newcommand{\clF}{\mathcal{F}}

\newcommand{\clH}{\mathcal{H}}
\newcommand{\clI}{\mathcal{I}}

\newcommand{\clP}{\mathcal{P}}

\newcommand{\clT}{\mathcal{T}}

\newcommand{\gC}{\mathfrak{C}}

\newcommand{\ga}{\mathfrak{a}}

\newcommand{\gd}{\mathfrak{d}}

\newcommand{\eps}{\varepsilon}
\newcommand{\cl}{c\kern-.11ex \ell}
\newcommand{\lh}{{\ell\kern-.27ex g}}
\newcommand{\rest}{\restriction}
\newcommand{\tser}{\upharpoonleft}

\newcommand{\caret}{{\char 94}}
\newcommand{\LL}{\langle}
\newcommand{\RR}{\rangle}

\newcommand{\subref}[1]{$_{\mathrm{\texttt{=}}\mathsf{L{#1}}}$}

\newcommand{\lepref}[1]{({<}\,{#1})}

\newcommand{\overbar}[1]{\mkern 1.5mu\overline{\mkern-1.5mu#1\mkern-1.5mu}\mkern 1.5mu}

\newcommand{\ol}{\overline}
\newcommand{\olsi}[1]{\,\overline{\!{#1}}} 

\newcommand*{\defeq}{\mathrel{\vcenter{\baselineskip0.5ex \lineskiplimit0pt\hbox{\scriptsize.}\hbox{\scriptsize.}}}=}

\usepackage{tcolorbox}

\newcount\skewfactor
\def\mathunderaccent#1#2 {\let\theaccent#1\skewfactor#2
\mathpalette\putaccentunder}
\def\putaccentunder#1#2{\oalign{$#1#2$\crcr\hidewidth
\vbox to.2ex{\hbox{$#1\skew\skewfactor\theaccent{}$}\vss}\hidewidth}}
\def\name{\mathunderaccent\tilde-3 }
\def\Name{\mathunderaccent\widetilde-3 }

\newbox\noforkbox \newdimen\forklinewidth
\setbox0\hbox{$\textstyle\bigcup$}
\setbox1\hbox to \wd0{\hfil\vrule width \forklinewidth depth \dp0
                        height \ht0 \hfil}
\setbox\noforkbox\hbox{\box1\box0\relax}
\def\unionstick{\mathop{\copy\noforkbox}\limits}
\def\nonfork#1#2_#3{#1\unionstick_{\textstyle #3}#2}
\def\nonforkin#1#2_#3^#4{#1\unionstick_{\textstyle #3}^{\textstyle
    #4}#2}
\setbox0\hbox{$\textstyle\bigcup$}
\setbox1\hbox to \wd0{\hfil{\sl /\/}\hfil}
\setbox2\hbox to \wd0{\hfil\vrule height \ht0 depth \dp0 width
                                \forklinewidth\hfil}
\newbox\doesforkbox
\setbox\doesforkbox\hbox{\box1\box0\relax}
\def\nunionstick{\mathop{\copy\doesforkbox}\limits}

\def\fork#1#2_#3{#1\nunionstick_{\textstyle #3}#2}
\def\forkin#1#2_#3^#4{#1\nunionstick_{\textstyle #3}^{\textstyle
    #4}#2}

\newcommand{\stickT}{%
\setbox255=\hbox{\raise1ex\hbox{$\hspace{0.2pt}\,\bullet\,$}}
\mathord{\rlap{\hbox to\wd255{\hss\hbox{$|$}\hss}}
\box255}
}
\newcommand{\stickS}{%
\setbox255=\hbox{\raise0.6ex\hbox{$\scriptstyle\bullet$}}
\mathord{\rlap{\hbox to\wd255{\hss\hbox{$\scriptstyle|$}\hss}}
\box255}
}

\author[Shelah]{Saharon Shelah}
\address{Einstein Institute of Mathematics,
The Hebrew University of Jerusalem,
9190401, Jerusalem, Israel; and\\
Department of Mathematics,
Rutgers University,
Piscataway, NJ 08854-8019, USA}
\urladdr{https://shelah.logic.at/}
\thanks{First typed 2022-03-25.
The author thanks Craig Falls  
for generously funding typing services, and Matt Grimes for the careful and beautiful typing.
The author would like to thank the Israel Science Foundation for partial support of
this research by grant 2320/23 (2023-2027).\\
References like (e.g.) [Sh:950, Th0.2\subref{y5}] mean that the internal label of Theorem 0.2 in Sh:950 is `\textsf{y5}.'
The reader should note that the version in my website is usually more up-to-date than the one in arXiv.
This is publication number
1257 on Saharon Shelah's list.
}

\makeatletter
\@namedef{subjclassname@2020}{\textup{2020} Mathematics Subject Classification}
\makeatother
\subjclass[2020]{Primary 03E35; Secondary 03E25, 03E15}
\keywords{set theory, forcing, iterated forcing, homogeneity, definability, axiom of choice, ZF+DC}
\date{March 10, 2026} 

\title[Homogeneous forcing]{Homogeneous forcing\\{}  1257}

\begin{document}
\makeatletter\def\shfiuwefootnote{\gdef\@thefnmark{}\@footnotetext}\makeatother\shfiuwefootnote{Version 2026-03-12. See \url{https://shelah.logic.at/papers/1257/} for possible updates.}
\begin{abstract}
Assume $\kappa = \kappa^{< \kappa}$ (usually $\aleph_0$ or an inaccessible).

We shall deal with iterated forcings preserving ${}^{\kappa>}{\rm Ord}$ and not collapsing cardinals along a linear order $L$. A sufficient condition for this, which we will focus on, is for the forcings to have support $<\kappa$ and the $\kappa^+$-cc, and be strategically $<{\kappa}$-complete. The aim is to have homogeneous forcings, so that the iteration has many automorphisms. 

In addition to the inherent interest, such iterations are helpful for considering some natural ideals on ${}^\kappa2$, in order to get a model of ${\rm ZF} + {\rm DC}_\kappa\ +$ ``modulo this ideal, every set is equivalent to a $\kappa$-Borel one."

But here we only have many automorphisms of the index set $L$ and therefore of the iteration of iterands $\mathbb{Q} $; 
we do not necessarily have homogeneity of $\mathbb{Q} $, and we do not have automorphisms mapping other names of $\mathbb{Q} $-reals onto each other.
However, for some reasonable forcing notions, there are no other $\mathbb{Q} $-reals! This was the reason for introducing and investigating saccharinity in earlier works with Jakob Kellner and with Haim Horowitz.
\end{abstract}

\maketitle

\centerline{\underline{\textbf{Annotated Content}}}

\sn
\textbf{\S0 \quad Introduction} \,$_{\mathsf{(label\ z)}}$ \hfill p.\pageref{S0}

\textbf{\S0}(A) \ Aim

\textbf{\S0}(B) \ Background

\textbf{\S0}(C) \ Preliminaries 

We define a trunk controller $\TR$ in \ref{z20}.

\bn
\textbf{\S1 \quad The Frame} \,$_{\mathsf{(label\ f)}}$ \hfill p.\pageref{S1}

\begin{quotation}
    We deal with the individual $\bbQ_\bfo$-s (for $\bfo \in \bfO_\TR$) and 
    $\bfm \in \bfM$. They will underpin the combinatorics of our abstract iteration.
\end{quotation}

\bn
\textbf{\S2 \quad Investigating $\bfp \in \bfP$} \,$_{\mathsf{(label\ f23{+})}}$ \hfill p.\pageref{S2}

\begin{quotation}
    We define $\bfP$, a class of abstract iterations (but {limits of} increasing sequences are well-defined only in the pure version).
\end{quotation}

\bn
\textbf{\S3 \quad From $\bfM$ and $\bfP$ to $\bfN$ and $\bfQ$} \,$_{\mathsf{(label\ d)}}$ \hfill p.\pageref{S3}

\begin{quotation}
    Now we have $W_\bfq \subseteq L_\bfq$ (for iterations $\bfq \in \bfQ$) such that for $t \in W_\bfq$, the memory $\clA_t = \clB_t$ behaves as in $\bfP$. (It even bounds the cardinality, and they are not changed for $\leq_\bfQ$-larger iterations.) 

    But for $t \in L_\bfq \setminus W_\bfq$, we have two kinds of memory $\clB_t \subseteq \clA_t$, where the $\clB_t$-s behave as above under the order on $L$ (and are used to show $\bfq_1 <_\bfQ \bfq_2 \Rightarrow \bbP_{\bfq_1} \lessdot \bbP_{\bfq_2}$).
\end{quotation}

\bn
\textbf{\S4 \quad Homogeneity} \,$_{\mathsf{(label\ g)}}$ \hfill p.\pageref{S4}

\bn
\textbf{\S5 \quad Variants, and the Theorem} \,$_{\mathsf{(label\ i)}}$ \hfill p.\pageref{S5}

\begin{quotation}
    The theorem says that saccharinity of $\bbQ_\bfo$ is sufficient to prove the consistency promised in \S0A-B.
\end{quotation}


\newpage
\setcounter{section}{-1}
\section{Introduction}\label{S0}

\subsection{Aim}

We fix $\kappa = \kappa^{< \kappa}$ and consider homogeneous $\lepref{\kappa}$-support iterations of 
$\lepref{\kappa}$-complete forcing notions, with a version of $\kappa^+$-cc, preserving those properties. However, throughout this section we will concentrate on the classical case $\kappa \defeq \aleph_0$.

To get homogeneity we intend to iterate along a linear order which is quite homogeneous (and therefore very much not well-ordered).

Ever since Solovay's celebrated work \cite{So70}, we know about the connection between the following two issues:
\begin{enumerate}[$\bullet_1$]
    \item  Forcing notions $\bbP$ with lots of automorphisms. E.g.~for small 
    $\bbP' \lessdot \bbP$ and two relevant $\bbP$-names $\name\eta_1,\name\eta_2$ of reals, generic for the same relevant forcing $\bbQ$ over 
    $\bfV^{\bbP'}$, there is an automorphism of $\bbP$ over $\bbP'$ mapping 
    $\name\eta_1$ to $\name\eta_2$.
\sn
    \item Models of $\ZF + \DC +$ ``every set of reals is equivalent to a Borel set modulo the null ideal (or other reasonable ideal)". (The most central forcings were Random Real forcings for the null ideal. The second prominent case {was} Cohen forcing, for the meagre ideal.)
\end{enumerate}
Concerning the classical case of Lebesgue measurability, another formulation is ``no non-measurable set is easily definable," formulated\footnote{
    That is, $\bullet_2$ holds for an inner model $\bfL[\clP(\kappa)]^\bfV$ with $\bfV \models \ZFC$, so in $\bfV$ all `reasonable' sets are `measurable' for this ideal.
}
as `definable in $\bfL[\bbR]$.' See the history and more in \cite{Sh:672}, \cite{Sh:856}.

This applies to other ideals $\id(\bbQ,\name\eta)$ for a definable forcing notion
$\bbQ$ (mainly a ccc one) and a $\bbQ$-name $\name\eta$ of a real. 
Generally, it
was not so easy to build such forcing notions: it required one to prove the existence of amalgamation in the relevant class of forcings. In Kellner-Shelah \cite{Sh:859} it was suggested to look at so-called saccharine pairs
$(\bbQ,\name\eta)$, where $\bbQ$ is very non-homogeneous. (E.g.~forcing with $\bbQ$ adds just one $(\bbQ,\name\eta)$-generic, so we have few cases we need to build automorphisms for. The forcing notion {here} is proper but \emph{not} ccc.)

Here we prove the existence of sufficiently homogeneous $\kappa^+$-cc 
$\lepref{\kappa}$-complete iterated forcings (so, along a linear order with $\lepref{\kappa}$-support). 

\mn
\begin{notation}\label{z2}
1) 
\begin{itemize}
    \item $\kappa,\lambda,\mu,\partial,\theta,\sigma$ are cardinals (infinite, if not explicitly said otherwise). $\lambda^+ = \lambda(+)$ will denote the successor of $\lambda$.
\sn
    \item $\alpha,\beta,\gamma,\delta,\eps,\zeta$ will denote ordinals; $\delta$ will be a limit ordinal if not stated otherwise.
\sn
    \item $k,\ell,m,n$ {were originally} natural numbers, but now just ordinals $< \kappa$.
\sn
    \item $i$ and $j$ have been used for indices $<\kappa$ throughout.
\end{itemize}

\sn
2) $S_\kappa^\lambda \defeq \{\delta < \lambda : \cf(\delta) = \kappa\}$.

\sn
3) Recall that $\bbL_{\kappa,\sigma}$ is defined like first-order logic, but allowing $\bigwedge\limits_{i< \alpha}\varphi_i$ for $\alpha < \kappa$ and $(\exists \ldots x_i\ldots)_{i\in I}$ with $I$ of cardinality $< \sigma$.

\sn
4) $\bbP,\bbQ,\bbR$ will denote forcing notions.
\end{notation}

\mn
\begin{definition}\label{z5}
0) For $\kappa$ a cardinal,
the family of $\kappa$-Borel sets is the smallest family of subsets of $2^\kappa$ containing all
basic sets of the form $\{\nu\in {}^\kappa2 : \nu(\alpha)=i\}$
and closed under complements and unions of ${\le}\,\kappa$-many sets.

\sn
1) For $\bbQ$ a forcing notion, $\name\eta$ a $\bbQ$-name of a member of $2^\kappa$, and $\partial$ a cardinal, let $\id_{< \partial}(\bbQ,\name\eta)$ be the ideal consisting of the unions 
of ${<} \, \partial$-many $\kappa$-Borel sets $\bfB$ such that 
$\Vdash_\bbQ ``\name\eta \notin \bfB"$.

\sn
2) We say $\id_{< \partial}(\bbQ,\name\eta)$ has \emph{measurability} when {for} every $Y \subseteq {}^\kappa2$ there exists a $\kappa$-Borel set $\bfB$ such that $Y\, \triangle\, \bfB \in \id_{< \partial}(\bbQ,\name\eta)$.

\sn
3) Let $\id_{\leq\partial}(\bbQ,\name\eta)$ be $\id_{< \partial^+}\!(\bbQ,\name\eta)$.
\end{definition}

\bn
\subsection{Background}\label{0b}

\begin{discussion}\label{z6}
1) Comparing the forcing $\bbQ$ from \cite{Sh:859} to the older results (such as Solovay \cite{So70}), the forcings are Borel definable, proper, and of cardinality $2^{\aleph_0}$. In addition:
\begin{enumerate}
    \item [$\bullet_{1.1}$] 
    The forcing $\bbQ$ collapsed no cardinal (provided that \textsf{CH} held), but was not ccc; this\footnote{
        Note that Solovay uses Levy collapse of an inaccessible, but the later versions use ccc ones (mainly for the meagre ideal).
    }
    we consider a drawback.
\sn
    \item [$\bullet_{1.2}$] The model, as in those older results, does satisfy $\ZF + \DC$.
\sn
    \item[$\bullet_{1.3}$] The iteration was along a homogeneous linear order.
\sn
    \item[$\bullet_{1.4}$] We get only a somewhat weaker version of measurability, the ideal being $\id_{\leq\aleph_1}\!(\bbQ,\name\eta)$ instead of
    $\id_{< \aleph_1}\!(\bbQ,\name\eta)$.
\end{enumerate}
Alternatively (e.g.~starting with the universe $\bfL$),
\begin{enumerate}
    \item[$\bullet_{1.4}'$] Use $\id_{< \aleph_1}\!(\bbQ,\name\eta) +X$, where $X$ is the set $\{\name\eta[\bfG] : \bfG \subseteq \bbQ^\bfL \text{ is generic over } \bfL\}$.
\end{enumerate}

\mn
2) The next step was Horowitz-Shelah \cite{Sh:1067}, where:
\begin{enumerate}
    \item[$\bullet_{2.1}$] The forcing is ccc, which we consider a plus.
\sn
    \item [$\bullet_{2.2}$] The model only satisfies $\ZF$; we do not get $\DC$ or even $\AC_{\aleph_0}$ --- not so good.
    \item[$\bullet_{2.3}$]  Again, the iteration is along a homogeneous linear order.
\sn
    \item[$\bullet_{2.4}$] The ideal is again
    $\id_{\leq\aleph_1}\!(\bbQ,\name\eta)$ (or as in $\bullet_{1.4}'$ above).
\end{enumerate}
\end{discussion}

\medskip
Our intention is to regain both ccc (as in $\bullet_{2.1}$) as well as $\DC$ (as in $\bullet_{1.2}$) for the ideal $\id_{\leq\aleph_1}(\bbQ,\name\eta)$. Moreover, we can demand $\DC_{\aleph_1}$ (or more; see \S1) which is a significant plus.

\mn
\centerline{*\qquad*\qquad*}

\medskip
We continue \cite{Sh:700}, \cite{Sh:700a}, but do not rely on them. Instead of defining iterations we introduce them axiomatically and
allow  $\kappa > \aleph_0$ (referenced in the completeness and the support).
Unlike \cite{Sh:700}, {the} present paper does not address forcing 
$\ga > \gd$. Earlier continuations of \cite{Sh:700} and \cite{Sh:700a} 
were the parallel papers, in preparation, with preliminary numbers F2009 and F2029 (and later, their descendants F2330 and F2329). 
In \cite{Sh:700} the set $\cP_{\!s}^\bfm$ (see Def.\ \ref{f11}) may be the 
whole power set, and we use more general definable forcing notions.

In our iteration we are allowed to replace $\aleph_0$ by some $\kappa = \kappa^{< \kappa}$, so the forcing notions are $\lepref{\kappa}$-complete $\kappa^+$-cc. For the $\kappa  =\aleph_0$ case we intend to use the forcing notion $\bbQ_\bfn^2$ from \cite{Sh:1067}. But we still need an analogue which works for any $\kappa$.

\bn
\centerline{\textbf{Explanation of the path chosen}}

(The reader may skip this subsection for now, as it will only make sense after reading \S1-2.)

Consider a family $\bfK$ of $\kappa^+$-cc $\lepref{\kappa}$-complete forcing notions, ordered by some $\leq_\bfK$. There are two natural ways to do this.

\mn
\underline{\textsc{The first way:}}

We build a directed family $\LL \bbP_t : t \in I\RR$ such that $s <_I t \Rightarrow \bbP_{\!s} \lessdot \bbP_t$ and $\bbP \defeq \bigcup\limits_{t\in I} \bbP_t$ is in $\bfK$ and has enough automorphisms.

Presently, 
this requires
\begin{enumerate}
    \item Start with $(\bfP,\leq_\bfP)$ with the entire partial order $<_\bfP^\pr$, as in the present version of \S2. (So no $\clB_\bfq$!)
\sn
    \item We define $\bfQ$ as the set of $\bfq$ which consist of:
    \begin{enumerate}
        \item $W_\bfq \subseteq L_\bfq$ such that for every $s \in L_\bfq$ we have $\clA_{\bfq,s}^\circ \subseteq \clA_{\bfq,s}$ with $\lambda_{\bfo_s}$ members, which are pairwise disjoint sets all of cardinality $\leq \lambda_{\bfo_s}$, and are \emph{dense} in some suitable sense.

        (Alternatively, allow $2^{\lambda_s}$-many such that all possible isomorphism types appear.)
\sn
        \item  For each $s \in L_\bfq \setminus W_\bfq$ we have an ideal $\clI_s$ on $\clA_{\bfq,s}^\circ$ such that 
        $$
        A \in \clA_{\bfq,s} \Rightarrow \{B \in \clA_{\bfq,s}^\circ : A \cap B \neq \varnothing\} \in \clI_s.
        $$
        \item If $s \in W_\bfq$ then $\bigcup\clA_{\bfq,s} \subseteq W_\bfq$.
\sn
        \item If $\bfq_1 \leq_\bfQ \bfq_2$ and $s \in W_{\bfq_1}$, \underline{then} $\clA_{\bfq_1,s} = \clA_{\bfq_2,s}$.
    \end{enumerate}
\end{enumerate}

\mn
\underline{\textsc{The wrong way:}}

Let $\chi = \chi^{<\chi} > \theta = \cf(\theta) > \kappa = \cf(\kappa)$ such that $\kappa = \kappa^{<\kappa}$ and $\alpha < \theta \Rightarrow |\alpha|^\kappa < \theta$. We force by $(\bfK_{<\theta} \cap \clH(\chi),<_\bfK)$, where $\bfK \defeq \bfP$ and $\bfK_{<\theta} \defeq \{\bbP \in \bfK : \|\bbP\| < \theta\}$.

To deal with $\bbL[\bbR^{<\theta}]$, it is enough to prove that there are enough amalgamation bases in $(\bfK_{<\theta},\leq_\bfK)$, which is not so hard. (It is enough to have ``the union of a $\leq_\bfK$-increasing chain of length $\partial$ is a lub," for some $\partial \in [\kappa,\theta) \cap \Reg$.)

A suspicious point is how to define when a (not overly long) increasing sequence in $(\bfP,\leq_\bfP)$ {has an upper bound,} even in the case $\kappa \defeq \aleph_0$. However, for $(\bfP,\leq_\bfP^M)$ this is not a problem.

Too good to be true? No; just that for homogeneity we need to return to the first way.

\mn
\centerline{*\qquad*\qquad*}

\medskip
We thank Shimon Garti, Miguel Cardona, Martin Goldstern, and Jakob Kellner for doing much to improve this paper.

\medskip

\subsection{Preliminaries}

\begin{hypothesis}\label{z8}
1) $\kappa = \kappa^{< \kappa}$ (mainly $\aleph_0$ or an inaccessible).

\sn
2) $\partial$ is a regular cardinal $> \kappa$.

\sn
3) $D$ is a normal filter on $\kappa^+$ such that $S_\kappa^{\kappa^+} \defeq \{\delta < \kappa^+ : \cf(\delta) = \kappa\} \in D$.
\end{hypothesis}

\bigskip
For more on the condition defined below, see \cite{Sh:1036}.
\begin{definition}\label{z11} 
1) For $D$ a normal filter on $\kappa^+$ containing $S_{\kappa}^{\kappa^+}$,
we say the forcing notion $\bbQ$ satisfies $*_{\kappa,D}^0$ \underline{when}:

\begin{enumerate}
    \item[$*_{\kappa,D}^0$] Given a sequence $\bar p = \LL p_i : i < \kappa^+\RR$ of members of 
    $\bbQ$, there is a set\footnote{
        Yes! Not just `$C \in D^+$;' see \cite{Sh:1036}.
    } 
    $C \in D$, a sequence $\bar p^+ = \LL p_i^+ : i< \kappa^+\RR$ satisfying $(\forall i) [p_i \le p_i^+]$, and a regressive function $\bfh$ on $C$ such that
    $$
    \alpha,\beta \in C \wedge \bfh(\alpha) = \bfh(\beta) \Rightarrow `p_\alpha^+ \text{ and $p_\beta^+$ have a lub.'}
    $$
\end{enumerate}

\mn
2) For $\kappa$ and $D$ as above, we say $\bbQ$ satisfies $*_{\kappa,D}^1$ \underline{when} $*_{\kappa,D}^0$ holds with $\bar p^+ \defeq \bar p$.


\mn
3) We say a forcing notion $\bbQ$ is \emph{Knaster}$^+$ \underline{when}:
\begin{quotation}
    If $p_i \in \bbQ$ for $i < \omega_1$ \underline{then} there are $A \in [\omega_1]^{\aleph_1}$ and $\LL p_i^+ : i \in A\RR \subseteq \bbQ$, with each $p_i^+ \geq_\bbQ p_i$, such that for every $i,j \in A$ the conditions $p_i^+$ and $p_j^+$ are compatible and have a least upper bound.
\end{quotation}
\end{definition}

\mn
\begin{notation}
1) $\bfp$ will be a member of $\bfP$ as in Definition \ref{f23}.

\sn
2) $\bfq,\bfr$ will denote ATIs (\emph{abstract template iterations}); i.e.~members of $\bfQ_\pre$ (the weakest version --- see Definitions \ref{d20}, \ref{g5}).

\sn
3) $L$ is a linear order (usually $L \subseteq L_\bfm$) and $r,s,t \in L$.

$L_+$ is derived from $L$, with $\infty$ and $t,t(+) \in L_+$ for each $t \in L$. (See below in \ref{f11}(2).)

\sn
4) $L_\bfm$ or $L_\bfq$ will be the relevant linear order for $\bfm$ or $\bfq$, etc.

\sn
5) $\bbP,\bbQ,\bbR$ denote forcing notions as in Definition \ref{z11} (which means quasi-orders).
\end{notation}

\mn
\begin{definition}\label{z20}
1) We say $\TR$ is a $\kappa$-\emph{trunk controller} \underline{when} it consists of:
\begin{enumerate}
    \item A partial order $\leq_\TR$ on a set $|\TR|$ of cardinality $\|\TR\|$. 
    
    (Abusing notation slightly, we may write `$s \in \TR$' instead of $s \in |\TR|$.)
\sn
    \item A partial function $\plus_\TR : {}^{\kappa>}|\TR| \to |\TR|$ giving a lub.
\sn
    \item Any increasing sequence of length $< \kappa$ belongs to $\dom(\plus_\TR)$ (and hence has a lub).
\sn
    \item If $\zeta < \kappa$ and $\bfx_{\alpha,\eps} \in \TR$ for $\alpha < \kappa^+$ and $\eps < \zeta$, \underline{then} for some $\alpha < \beta$ we have 
    $$
    \bigwedge_{\eps< \zeta} \big[ \LL\bfx_{\alpha,\eps},\bfx_{\beta,\eps}\RR \in \dom(\plus_\TR) \big].
    $$
    \item $\varnothing \in \TR$, and $(\forall \bfx \in \TR) \big[ \varnothing \leq_\TR \bfx\big]$.
\sn
    \item $S_\TR \subseteq \clH_{<\kappa}(\Ord)$ and 
    $$\val_\TR : |TR| \times S_\TR \to 2$$ 
    is a partial function. 
    
    Specifically, for each $\bfx \in |\TR|$, $\val(\bfx) = \val(\bfx,-)$ is a partial function from $S_\TR$ to $\{0,1\}$, increasing with $\bfx$.
\end{enumerate}

\sn
1A) 
We call $\TR$ a \emph{simple} $\kappa$-trunk controller \underline{when}:
\begin{enumerate}
    \item $\TR$ is a $\kappa$-trunk controller.
\sn
    \item $(\forall \bfx \in \TR)\big[ |\bfx| \leq \kappa \big]$ (We usually have  $|\TR|\subseteq \clH(\kappa)$, or at least $\clH_{<\kappa}(\lambda)$ for some $\lambda$.)
\sn
    \item ${\leq_\TR} \defeq {\subseteq}$
\sn
    \item If $\delta < \kappa$ is limit and $\LL \bfx_\alpha : \alpha < \delta\RR$ is $\leq_\TR$-increasing, \underline{then}
    $$
    \plus_\TR(\LL \bfx_\alpha : \alpha < \delta\RR) \defeq \bigcup\limits_{\alpha<\delta} \bfx_\alpha.
    $$
    \item $S_\TR \defeq |\TR|$, and for each $\bfx \in \TR$, we define the function $\val_\TR(\bfx)$ by demanding that $\val_\TR(\bfx)(\bfy)$ is 1 if $\bfy \leq_\TR \bfx$ and is 0 if $\bfx$ and $\bfy$ are incompatible (and is undefined otherwise). 
\end{enumerate}

\sn
2) We say $\TR$ is a $(\kappa,\Theta,\Upsilon)$-{trunk controller} \underline{when}:
\begin{enumerate}[$\bullet_1$]
    \item $\Theta$ is a set of regular cardinals $<\kappa$.
\sn
    \item $\kappa >\aleph_0 \Leftrightarrow \Theta \neq \varnothing$
\sn
    \item $\Upsilon \in [2,\kappa]$.
\sn
    \item If $\Theta = \Reg \cap \kappa$ then we may omit it; similarly if $\Upsilon = 2$.
\end{enumerate}
\begin{enumerate}
    \item [(A),] (B), {(C),} (E) As above.
\sn
    \item [(D)$'$] Any $\leq_\TR$-increasing sequence $\LL \bfx_i : i < j\RR$ with $j < \kappa$ has an upper bound, \underline{but} if $j \in \Theta$ then it has a lub $\plus_\TR(\LL \bfx_i : i < j\RR)$. 
\sn
    \item[(F)] If $j < 1+\Upsilon$, $\zeta < \kappa$, and $\bfx_{\alpha,\eps} \in \TR$ for $\alpha < \kappa^+$ and $\eps < \zeta$, \underline{then} for some increasing sequence of ordinals $\LL\alpha_i : i < j\RR$ bounded above by $\kappa^+$, for each $\eps < \zeta$, we have
    $$
    \LL \bfx_{\alpha_i,\eps} : i < j\RR \in \dom(\plus_\TR).
    $$
\end{enumerate}


\mn
3) For $\LL \TR_s : s \in I\RR$, with each $\TR_s$ as in (1A) or (2),
we define $\TR = \prod\limits_{s \in I}^{<\kappa} \TR_s$, the $\lepref{\kappa}$-\emph{support product}, as follows.
\begin{enumerate}
    \item $f \in \TR$ \underline{iff} $f \in \prod\limits_{s\in J} \TR_s$ for some $J \in [I]^{<\kappa}$.
\sn
    \item $f \leq_\TR g$ \underline{iff} ($f,g \in \TR$ and)
    \begin{enumerate}
        \item  $\dom(f) \subseteq \dom(g)$
\sn
        \item $f(s) \leq_{\TR_s} g(s)$ for all $s \in \dom(f)$.
    \end{enumerate}
\sn
    \item Let $\bar f = \LL f_\eps : \eps < \zeta\RR$, with $\zeta < \kappa$. 

\sn
    We have $\plus_\TR(\bar f) = g$ \underline{when}
    \begin{enumerate}
        \item $f_\eps \in \TR$ for all $\eps < \zeta$.
\sn
        \item $\dom(g) = \bigcup\limits_{\eps < \zeta} \dom(f_\eps)$ 
\sn
        \item For all $s \in \dom(g)$, if $u = u_{\bar f,s} \defeq \big\{ \eps < \zeta : s \in \dom(f_\eps)\big\}$, \underline{then} 
        $$
        g(s) = \plus_{\TR_s}(\LL f_\eps(s) : \eps \in u\RR).
        $$
    \end{enumerate}
\sn
    \item 
    \begin{enumerate}
        \item $S_\TR \defeq \big\{ (s,\bfx) : s \in I,\ \bfx \in S_{\TR_s}\big\}$
\sn
        \item For $f \in \TR$, define $\val(f)$ as the function with domain 
        $$
        \big\{ (s,\bfx) : s \in \dom(f),\ \bfx \in \dom(\val(f(s)))\big\}
        $$ 
        which sends $(s,\bfx) \mapsto \val_{\TR_s}(f(\bfx))$. 
    \end{enumerate}
\end{enumerate}
\end{definition}

\mn
\begin{remark}\label{z25}
1) For $\TR$ a $(\kappa,\Theta,\Upsilon)$-trunk controller, the function $\plus_\TR$ will not really be used for sequences of length $j \in \kappa \cap \Reg \setminus \Theta$; only the existence {of such a value will be necessary.} (See more on this in \S5.)

\mn
2) If $\TR_s$ is a $\kappa$-trunk controller for $s \in I$, then 
$\prod\limits_{s \in I}^{<\kappa} \TR_s$ is not \emph{guaranteed} to be a $\kappa$-trunk controller, as \ref{z20}(1)(D) may fail.
\end{remark}

\mn
\begin{example}\label{z28}
We define an explicit $\kappa$-trunk controller $\TR = \TR_\kappa$ as follows.
\begin{enumerate}
    \item The set of elements is $\{f \in \clH(\kappa) : f \text{ is a function}\}$.
\sn
    \item ${\leq_\TR} \defeq {\subseteq}$
\sn
    \item $\plus_\TR(\LL f_i : i < j\RR) \defeq g$ \ \underline{iff} \ $\bigwedge\limits_{i<j} [f_i = g \in \TR]$.
\sn
    \item $S_\TR \defeq |\TR|$ and 
$$
    \dom(\val_\TR(f)) \defeq \big\{g \in |\TR| : \dom(g) \subseteq \dom(f)\big\}
$$
    with 
$$
    \val_\TR(f)(g) \defeq
    \begin{cases}
        1 &\text{if } g \subseteq f\\
        0 &\text{otherwise.}
    \end{cases}
$$
\end{enumerate}
\end{example}

\mn
\begin{claim}\label{z31}
$1)$ Assume $\kappa = \kappa^{<\kappa}$. If $I$ is a set and $\TR_s$ 
is a $\kappa$-trunk controller for each $s \in I$, \underline{then} 
$\TR \defeq \prod\limits_{s \in I}^{<\kappa} \TR_s$ is a $\kappa$-trunk controller \underline{iff} $\prod\limits_{s \in J}^{<\kappa} \TR_s$ is a $\kappa$-trunk controller for all $J \in [I]^{<\kappa}$.

\mn
$2)$ Above, if $\TR_s$ is constant in $s$ \underline{then} the product is also a $\kappa$-trunk controller.

\mn
$3)$ $\TR_\kappa$ from Example \emph{\ref{z28}} is indeed a $\kappa$-trunk controller.

\mn
$4)$ Every $\kappa$-trunk controller is a $(\kappa,\Reg \cap \kappa,2)$-trunk controller.
\end{claim}

\begin{PROOF}{\ref{z31}}
Straightforward.
\end{PROOF}

\newpage

\section{The frame}\label{S1}

To explain our framework, we shall first concentrate on the most simple version: {every member of} $\bfTR$ is simple (see \ref{z20}(1A)) and each 
$\clA_t$ is a family of pairwise disjoint sets.



\begin{hypothesis}\label{f2}
For the rest of this work, we shall assume the following:

\sn
1) $\kappa = \kappa^{<\kappa}$

\sn
2)  $\bfTR$ is a non-empty set (or class) of {definitions of} $\kappa$-trunk controllers, such that the product of any $<\kappa$ of them is a $\kappa$-trunk controller.

\sn
3) $\bfTR$ is \emph{simple}; that is, every member of $\bfTR$ is a simple 
$\kappa$-trunk controller (see \ref{z20}(1A)).  

But we may use $S_\TR$ (or $S_\bfo$) and 
$\val_\bfo$ below  (see \ref{f32}).

\sn
4) $\bfTR^+ = \bfTR(+)$ is the closure of $\bfTR$ under products of length $< \kappa$.

\sn
5) If $(\forall s \in I)[\TR_s \in \bfTR]$ \underline{then} $\prod\limits_{s \in I}^{<\kappa} \TR_s$ is a $\kappa$-trunk controller. (Note that this is an easy consequence of part (2) above.)

\sn
6) $\bfO$ is a subset (or subclass) of $\bfO_\bfTR$ --- see Definition \ref{f3}. 
\end{hypothesis}

\mn
\begin{definition}\label{f3}
$\bfO_\bfTR$ is the class of objects $\bfo$ 
which consist of:\footnote{
    In the future, we may demand that {$\bfo$ has} a sufficiently absolute definition from the parameter $v \in {}^{W(\bfo)}2$ --- but for now, $W(\bfo) \defeq \varnothing$ is enough.
}
\begin{enumerate}
    \item $\TR_\bfo = \TR[\bfo] = (|\TR_\bfo|,\leq_\bfo,\plus_\bfo)$ is a $\kappa$-trunk controller from $\bfTR$ (and hence is simple, \emph{per} our hypothesis).
\sn    
    \item The \emph{atomic forcing} $\bbQ_\bfo^\at = (Q_\bfo^\at,\leq_\bfo^\at)$, which is a definition of a forcing notion (i.e.~a quasiorder).
\sn    
    \item 
    \begin{enumerate}
        \item A function $\tr_\bfo : Q_\bfo^\at \to |\TR_\bfo|$.
\sn
        \item  A non-decreasing function $\sigma_\bfo : |\TR_\bfo| \to [2,\kappa]_\Card$.
\sn
        \item  A function $A = A_\bfo$ with domain $\bbQ_\bfo^\at$ and range $\subseteq \clP(\TR_\bfo)$. We may write $A_{\bfo,p}$ in place of $A_\bfo(p)$.
    \end{enumerate}    
    (So again, these are \emph{definitions of functions}.)
\sn
    \item \textbf{[Notation:]} 
    \begin{enumerate}
        \item Let $\sigma_\bfo(p) \defeq \sigma_\bfo(\tr_\bfo(p))$ for $p \in Q_\bfo^\at$.
\sn
        \item $\lambda_\bfo \defeq \big( |Q_\bfo^\at| + \|\TR_\bfo\| + \kappa\big)^\kappa$
\sn
        \item For $\bfO \subseteq \bfO_\bfTR$, let $\lambda_\bfO \defeq \big( \sum\limits_{\bfo \in \bfO} \lambda_\bfo \big)^\kappa$.
    \end{enumerate}
\end{enumerate}
In addition, we make the following demands that the constituents of $\bfo$ must satisfy.
\begin{enumerate}
    \item[(E)]
    \begin{enumerate}
        \item Each $p \in Q_\bfo^\at$ is a subset of $\clH(\kappa)$ (or just $\subseteq \clH_{<\kappa}(\lambda_\bfo)$).
\sn
        \item If $p \leq_\bfo q$ \underline{then} $\tr(p) \leq_{\TR[\bfo]} \tr(q)$.
    \end{enumerate}
\sn 
    \item[(F)]
    \begin{enumerate}
        \item Let us say that a sequence $\bar p$ is $\sigma_\bfo$-\emph{truly increasing} \underline{when}
        $$
            i < j < \lh(\bar p) \Rightarrow \big[\sigma_\bfo(p_i) < \sigma_\bfo(p_j) \vee \sigma_\bfo(p_i) = \sigma_\bfo(p_j) = \kappa\big].
        $$
        \item If $\delta < \kappa$ is a limit ordinal and $\bar p = \LL p_\eps : \eps < \delta\RR$ is $\sigma_\bfo$-{truly increasing}, 
        \underline{then} $\bar p$ has a $\leq_\bfo$-lub
$$
        \textstyle\bigcup \bar p \defeq \big( \plus_\bfo(\LL\tr(p_\eps) : \eps < \delta \RR), \bigcup\limits_{\eps< \delta} p_\eps \big).
$$
        \item $A_{\bfo,p} \subseteq \{\tr(q) : q \ge_\bfo p\}$
    \end{enumerate}
\sn 
    \item[(G)] \textsc{Case 1:} [$\TR_\bfo \subseteq \clH(\kappa)$ for all $\bfo \in \bfO$.]

    If $\bfx \in \TR_\bfo$ and $\tr_\bfo(p_i) = \bfx$ for all $i < i_* < \sigma_\bfo(\bfx)$, \underline{then} $\{p_i : i < i_*\}$ has a common upper bound $q$ with $\tr(q) \in \bigcap\limits_{i<i_*} E_{\bfo,p_i}$.


    \mn 
    \textsc{Case 2:} [Holds in general.]

    If $\bfx_i \in \TR_\bfo$, $\LL\sigma_\bfo(\bfx_i) : i < i_*\RR$ is constant, $\tr_\bfo(p_i) = \bfx_i$ for all $i < i_* < \sigma_\bfo(\bfx_0)$, and $\LL\bfx_i : i < i_*\RR \in \dom(\plus_\bfo)$, \underline{then} $\{p_i : i < i_*\}$ has a common upper bound $q$ (again with $\tr(q) \in \bigcap\limits_{i<i_*} E_{\bfo,p_i}$).

\sn
    \item [(H)] \textbf{[Absoluteness:]} If $\bbR_1 \lessdot \bbR_2$ such that $\bbR_1$ and $\bbR_2/\bbR_1$ are strategically $\lepref{\kappa}$-complete $\kappa^+$-cc forcings, \underline{then}
    \begin{enumerate}
        \item In $\bfV^{\bbR_2}$, the definitions in $\bfo$ satisfy clauses (A)-(G). 
\sn
        \item  In $\bfV^{\bbR_2}$ we have $\bbQ_\bfo^{\bfV[\bbR_1]} \subseteq \bbQ_\bfo^{\bfV[\bbR_2]}$, 
        $\tr_\bfo^{\bfV[\bbR_1]} = \tr_\bfo^{\bfV[\bbR_2]}$, 
        $\sigma_\bfo^{\bfV[\bbR_1]} \subseteq \sigma_\bfo^{\bfV[\bbR_2]}$, and $A_\bfo^{\bfV[\bbR_1]} \subseteq A_\bfo^{\bfV[\bbR_2]}$.
\sn
        \item If $p_1 \in \bbQ_\bfo^{\bfV[\bbR_1]}$, {$p_2 \in \bbQ_\bfo^{\bfV[\bbR_2]}$,} $\bbQ_\bfo^{\bfV[\bbR_2]} \models `p_1\leq p_2$', and $\bfy \defeq \tr(p_2)$,
        \underline{then} there is $p_3 \in \bbQ_\bfo^{\bfV[\bbR_1]}$ such that $\tr(p_3) = \bfy$ and $\bbQ_\bfo^{\bfV[\bbR_1]} \models `p_1\leq p_3$'.
    \end{enumerate}
\sn
    \item $\name\eta_\bfo$ is the $\bbQ_\bfo^\at$-name of the function 
    $$
    \big\{ (s,\iota) \in S_\bfo \times 2 : (\exists p \in \Name\bfG_{\bbQ_\bfo})[\val(\tr(p),s) = \iota]\big\}.
    $$
\end{enumerate}
\end{definition}

\mn
The following definition will make more sense when used in our iterations.
\begin{definition}\label{f4}
For $\bfo \in \bfO_\bfTR$, we define $\bbQ_\bfo = (Q_\bfo,\leq_\bfo)$ as follows.
\begin{enumerate}
    \item We let $\tr_\bfo(p) \defeq \tr(p)$ and $\sigma_\bfo(p) \defeq \sigma_\bfo(\tr_p)$.
\sn
    \item Each $p \in Q_\bfo$ is of the form $(\tr_p,\clF_p)$, where $\tr_p = \tr(p) \in \TR_\bfo$ and $\clF_p$ is a set of members of $Q_\bfo^\at$ of cardinality $< \sigma_\bfo(\tr_p)$ such that (recalling \ref{f3}(F)(c))
    $$
    q \in \clF_p \Rightarrow \tr_\bfo(q) \leq_{\TR[\bfo]} \tr_\bfo(p) \in {A_{\bfo,q}}.
    $$
    \item For $p,q \in Q_\bfo$, we say $p \leq_{\bfo} q$ \underline{iff} $\tr(p) \leq_{\TR} \tr(q)$ and\footnote{
        Why not just ``for every $r_1 \in \clF_p$ there exists $r_2 \in \clF_q$ such that $r_1 \leq_\bfo r_2$?" Then we would only get $\lepref{\kappa}$-strategic completeness.
    } 
    $\clF_{p} \subseteq \clF_{q}$.
\sn
    \item When we say `$\LL p_i : i < \delta\RR$ is $\leq_\bfo$-\emph{truly increasing} in $\bbQ_\bfo$,' we mean it is increasing in $\bbQ_\bfo$ and 
    $\sigma_\bfo$-truly increasing as defined in \ref{f3}(F)(a).
\sn
    \item We define $\name\eta_\bfo$ as in Definition \ref{f3}. (Abusing our notation slightly, we will not distinguish between the two uses.)
\end{enumerate}
\end{definition}

\mn
\begin{claim}\label{f5}
Assume $\bfo \in \bfO_\bfTR$.

\mn
$1)$ $\bbQ_\bfo$ is a strategically $\lepref{\kappa}$-complete $\kappa^+$-cc forcing. 

\mn
$2)$ Assume $(\forall \bfo \in \bfTR)\big[\TR_\bfo \subseteq \clH(\kappa) \big]$ (\emph{or} just that for every $\mu < \kappa$, the set\\ $\{s \in \TR_\bfo : \sigma_\bfo(s) \geq \mu\}$ is {equal to} $\bigcup\limits_{\eps < \kappa} Y_\eps$ such that\footnote{
    A typical case is when $\tr(p)$ is a function from some $u \in [2^\kappa]^{<\kappa}$ to $\clH(\kappa)$ and $s \in \TR_\bfo \Rightarrow \sigma_\bfo(s) = \min\!\big\{\sigma_\bfo(s \rest \{\alpha\}) : \alpha \in \dom(s) \big\}$.
}
$\bar t \in {}^{\mu>}Y_\eps \Rightarrow \bar t \in \dom(\plus_\bfo)$).

For every $\mu < \kappa$, there exists a dense open subset of $\bbQ_\bfo$ which is the union of $\kappa$-many sets such that any $<\mu$ members of any one set has a common upper bound.

\mn
$3)$ If $\delta < \kappa$ is a limit ordinal and 
$\bar p = \LL p_i : i < \delta\RR$ is truly $\leq_\bfo$-increasing in 
$\bbQ_\bfo$, then $\bar p$ has a $\leq_{\bbQ_\bfo}$-lub (call it $p_\delta$) such that 
\begin{itemize}
    \item $\tr(p_\delta) = \plus_{\TR_\bfo}(\LL \tr_\bfo(p_i) : i < \delta\RR)$
\sn
    \item $\clF_{p_\delta} = \bigcup\limits_{i < \delta} \clF_{p_i}$.
\end{itemize}

\mn
$4)$ If $p_i \in Q_\bfo$, $\tr_\bfo(p_i) = \bfx_i$ for $i < i_*$, 
$\LL\bfx_i : i < i_*\RR \in \dom(\plus_\bfo)$, and $i_* < \sigma_\bfo(\bfx_0)$, \underline{then} there is a $\bfy \in \TR_\bfo$ such that 
$\big( \bfy, \bigcup\limits_{i< i_*} \clF_{p_i} \big) \in Q_\bfo$ is a common upper bound of $\{p_i : i < i_*\}$ in $\bbQ_\bfo$.

\mn
$5)$ The set $\{ p \in Q_\bfo : |\clF_p| \leq 1\}$ is dense in $\bbQ_\bfo$.
\end{claim}

\begin{PROOF}{\ref{f5}}
Straightforward. E.g.\ for part (1),

\mn
1) \textbf{[$\lepref{\kappa}$-complete:]}  Obvious.

\mn
\textbf{[$\kappa^+$-cc:]}

For each $p \in \bbQ_\bfo$, choose $\nu = \nu_p \in \TR$ such that $\tr(p) \leq_{\TR_\bfo} \bfx_p$ and $\sigma_\bfo(\bfx_p) > \mu \cdot |\clF_p|$. Suppose 
$\LL p_\alpha : \alpha < \kappa^+\RR$ is a sequence of members of $\bbQ_\bfo$ and 
$\bfx_\alpha \defeq \bfx_{p_\alpha}$ for $\alpha < \kappa^+$.

By the demand on $\TR_\bfo$, there are $\alpha_0 < \ldots < \alpha_{n-1}$ (for $n \geq 2$) such that $\plus_\bfo(\LL \bfx_{\alpha_\ell} : \ell < n\RR)$ 
is well defined
so easily $\{ p_{\alpha_\ell} : \ell < n\}$ has a common upper bound $q$ such that $\tr(q)$ is a lub of $\LL \bfx_{\alpha_\ell} : \ell < n\RR$.
\end{PROOF}

\mn
\begin{claim}\label{f7}
$1)$ If $\bfo \in \bfO$ then 
any $\sigma_\bfo$-truly increasing sequence of length $< \kappa$ has a lub. (See \emph{\ref{f3}(G)(c),(d)}.)

\mn
$2)$ If $\bfo_s \in \bfO$ (or just $\bfo_s \in \bfO^+$) for $s \in I$, \underline{then} $\bfo \defeq \prod\limits_{s \in I}^{<\kappa} \bfo_s$, naturally defined, belongs to $\bfO_{\bfTR(+)}$, and $\TR_\bfo = \prod\limits_{s \in I}^{<\kappa} \TR_s$.

\mn
$3)$ Assume $\bfo \in \bfO$.
\begin{enumerate}[$(A)$]
    \item $(\bfx,\clF) \in \bbQ_\bfo$ \underline{iff}
    \begin{enumerate}[$(a)$]
        \item $\bfx \in \TR_\bfo$
\sn
        \item $|\clF| < \sigma_\bfo(\bfx)$
\sn
        \item $f \in \clF \Rightarrow \big( \bfx,\{f\} \big) \in \bbQ_\bfo^\at$.
    \end{enumerate}
\sn
    \item $(\bfx_1,\clF_1) \leq_\bfo (\bfx_2,\clF_2)$ \underline{iff}
    \begin{enumerate}[$(a)$]
        \item $(\bfx_1,\clF_1), (\bfx_2,\clF_2) \in \bbQ_\bfo$
\sn
        \item $\clF_1 \subseteq \clF_2$
\sn
        \item $f \in \clF_1 \Rightarrow \big(\bfx_1,\{f\} \big) \leq_\bfo \big( \bfx_2,\{f\} \big)$.
    \end{enumerate}
\end{enumerate}
\end{claim}

\begin{PROOF}{\ref{f7}}
Easy.
\end{PROOF}

\mn
\begin{example}\label{f9}
1) Let $\bfn$ and $\bbQ_\bfn^2$ be as in \cite[Def.\,2.5]{Sh:1067}.
\begin{enumerate}
    \item We define an $\aleph_0$-trunk controller $\TR = \TR_\bfn = \TR[\bfn]$ as follows.
    \begin{enumerate}
        \item  $\kappa \defeq \aleph_0$
\sn
        \item $\TR_\bfn$ has set of elements $T_\bfn$ (as in \cite{Sh:1067}).
\sn
        \item $\sigma_{\TR_\bfn}$ will be the function $\eta \mapsto \mu_\eta$ (again, see there).
    \end{enumerate}
\sn
    \item We define $\bfO_\bfn \defeq \{\bfo : \bbQ_\bfo = \bbQ_\bfn^2\}$ (so it is a singleton), and $\bfO_{1067}$ will be the union of all $\bfO_\bfn$-s for $\bfn$ as in \cite[Def.\,2.5]{Sh:1067}.
\end{enumerate}

\mn
2) Let
$\kappa$ be inaccessible and $\bbQ_{\bar\lambda}$ be defined similarly to that in \cite[0.2\subref{z23}]{Sh:1126} (but there $\lambda_\eta \defeq \theta_{\lh(\eta)}$). It will be  of the form $\bbQ_\bfo$, where $\bfo = \bfo_{\bar\lambda} \in \bfO_\bfTR$. 

\noindent
More fully,
\begin{enumerate}[$\bullet_1$]
    \item $\bar\lambda = \LL\lambda_\eta, \lambda_\eta^0, D_\eta : \eta \in \clT\RR$
\sn
    \item $\clT \defeq \{\eta \in {}^{\kappa>}\!\kappa : \eps < \lh(\eta) \Rightarrow \eta(\eps) < \lambda_{\eta\rest\eps}\}$
\sn
    \item If $\eta \neq \nu \in \clT$ then $\lambda_\eta \neq \lambda_\nu$.
\sn
    \item If $\lh(\eta) < \lh(\nu)$ then $\lambda_\eta < \lambda_\eta^0  \leq \lambda_\nu < \kappa$.
\sn
    \item If $\eta \in \clT$ then $\cf(\lambda_\eta) \ge \lambda_\eta^0 = \cf(\lambda_\eta^0) > \prod\limits_{\lambda_\nu < \lambda_\eta} \lambda_\nu$.
\sn
    \item $D_\eta$ is a $\lambda_\eta^0$-complete uniform filter on 
    $\lambda_\eta$\\ (and naturally $A \in [\lambda_\eta]^{<\lambda_\eta} \Rightarrow \lambda_\eta \setminus A \in D_\eta$).
\end{enumerate}
Now let us work towards defining $\bfo = \bfo_{\bar\lambda}$.
\begin{enumerate}
    \item $\TR = \TR_{\bar\lambda} = \TR[\bar\lambda]$ is defined as follows.
    \begin{enumerate}
        \item $\kappa_\TR \defeq \kappa$
\sn
        \item $|\TR| \defeq \big\{\eta \in {}^{\kappa>}\!\kappa : \eps < \lh(\eta) \Rightarrow \eta(\eps) < \lambda_{\eta\rest\eps} \big\}$
\sn
        \item $\sigma_\TR(\eta) \defeq \lambda_\eta^0$.
    \end{enumerate}
\sn
    \item $p \in \bbQ_{\bar\lambda}^\at$ \underline{iff}
    \begin{enumerate}
        \item $p$ is a $\lepref{\kappa}$-complete subtree of $\TR$.
\sn
        \item $\tr(p)$ is the trunk of $p$.
\sn
        \item If $\tr(p) \unlhd \eta \in p$, then $\{\alpha < \lambda_\eta : \eta \caret \LL\alpha\RR \in p\} \in D_\eta$.
    \end{enumerate}
\sn
    \item $p \leq_{\bbQ_{\bar\lambda}} q$ \underline{iff} $p \supseteq q$.
\end{enumerate}

\sn
3) Assume $\kappa$ is inaccessible, and suppose $\clT = \clT_{\bar\lambda}$ is a subtree of $({}^{\kappa>}\!\kappa,\lhd)$ (so it is closed under initial segments and includes $\LL\ \RR$). Further suppose $\clT$ has no maximal node, and is closed under unions of $\lhd$-increasing chains of length $<\kappa$.

Let $\bar\lambda = \LL \lambda_\eta, \lambda_\eta^0, D_\eta : \eta \in \clT\RR$ be such that `$\eta \in \clT \Rightarrow (A)_\eta \vee (B)_\eta$' holds, where
\begin{enumerate}
    \item[$(A)_\eta$] $\lambda_\eta = \lambda_\eta^0 = 1$ and $D_\eta = \{0\}$.
\sn
    \item[$(B)_\eta$] $\lambda_\eta^0$ is an infinite regular cardinal $\leq \lambda_\eta < \kappa$, and $D_\eta$ is a $\lambda_\eta^0$-complete filter on $\lambda_\eta$.
\end{enumerate}
Furthermore, we demand that clause (2)\,$\bullet_2$ holds. 

Now define $\bfo_{\bar\lambda}$ as in {clause (2)}.

We say that $\bar\lambda$ (or $\clT_{\bar\lambda}$, or $\bfo_{\bar\lambda}$) is a $\kappa$-\emph{candidate} \underline{when} it is as above.

\sn
4) We say that $\bar\lambda$ (and $\bfo_{\bar\lambda}$) is $\kappa$-\emph{active} \underline{when}, in addition to the above, we have the following:
\begin{itemize}
    \item For every $\eta \in \lim(\clT)$, the set 
    $$
    \big\{ \delta < \kappa : \lambda_{\eta \rest\delta} > 1 \text{ and } (\forall \nu \in \clT \cap {}^\delta\kappa)[\nu \neq \eta \rest \delta \Rightarrow \lambda_{\eta\rest\delta} < \lambda_\nu \vee \lambda_\nu = 1 ] \big\}
    $$
    is stationary.
\sn
    \item If $\eta,\nu \in \clT$ with $\lh(\eta) < \lh(\nu)$, then 
    $\lambda_\nu \neq 1 \Rightarrow \lambda_\eta <\lambda_\nu^0$.
\end{itemize}
\end{example}

\mn
\begin{observation}\label{f10}
If $\kappa$ is inaccessible (or weakly inaccessible) and $\diamondsuit_\kappa$ holds, then there is a $\kappa$-active $\bar\lambda$.
\end{observation}

\mn
\begin{remark}
We will return to this example in \S3.
\end{remark}

\mn
\begin{definition}\label{f11}
0) Let $\bfM$ be the class of combinatorial templates (defined below). 

\sn
1) A \emph{combinatorial template} (or CT) $\bfm$ consists of:
\begin{enumerate}[(a)]
    \item A linear order $L_\bfm$ (we could have used `partial order'; it does not really matter for our purposes).

\sn
    We may write $s \in \bfm$ instead of $s \in L_\bfm$, or $s <_\bfm t$ instead of $s <_{L_\bfm} t$.
\sn
    \item  
    \begin{enumerate}[$\bullet_1$]
        \item $\olsi\clA = \LL \clA_t : t \in L_\bfm\RR$ (called the \emph{memory}), where $\clA_t \subseteq \clP(L_{\bfm,t})$ for each $t \in L_\bfm$ and $L_{\bfm,t} \defeq \{s \in L_\bfm : s <_{L_\bfm} t\}$.
\sn
        \item We denote $\clA \defeq \bigcup\limits_{t\in L_\bfm} \clA_t$.
    \end{enumerate}
\sn
    \item  {Each} $\clA_t$ is a family of subsets of $L_{\bfm,t}$ 
    with $\varnothing \in \clA_t$. We let $\clA_t^+ \defeq \clA_t \setminus \{\varnothing\}$.
    
\end{enumerate}


\sn
2) For $\bfm \in \bfM$, we {add} new objects $t(+)$ for all $t \in L_\bfm$, as well as $\infty$, and define $L_\bfm^+$, $L_{\bfm,x}$, $L_{\bfm,x}^+$, etc. as follows.
\begin{enumerate}[(a)]
    \item $L_\bfm^+ \defeq \{t,t(+) : t \in L_\bfm\} \cup \{\infty\}$ 
\sn    
    \item Naturally, $\LL t : t \in L_\bfm\RR \caret \LL t(+) : t \in L_\bfm\RR \caret \LL\infty\RR$ is without repetition.
\sn
    \item $<_{L_\bfm^+}$ is the closure, to a linear order, of the set
    $$
    \big\{ t < t(+) : t \in L_\bfm \big\} \cup \big\{ s(+) < t : s <_{L_\bfm} t \big\} \cup \big\{ t(+) < \infty : t \in L_\bfm \big\}.
    $$
    \item For $t \in L_\bfm^+$, let 
    $L_{\bfm,t} \defeq \{s \in L_\bfm : s <_{L_\bfm^+} t\}$ and 
    $$
    L_{\bfm,t}^+ \defeq \{s \in L_\bfm^+ : s <_{L_\bfm^+} t\}.
    $$
\end{enumerate}


\sn
3) For $L \subseteq L_\bfm$, we define $\bfm \rest L \in \bfM$ as follows.
\begin{enumerate}[$\bullet_1$]
    \item $L_{\bfm \rest L} \defeq L$
\sn
    \item $\clA_t^{\bfm \rest L} \defeq \clA_t^\bfm \cap \clP(L)$.
\end{enumerate}

\sn
4) For $t \in L_\bfm^+$, let $\bfm \rest t \defeq \bfm \rest L_{\bfm,t}^+$.

\sn
5) We say $\pi$ is an \emph{isomorphism from $\bfm_1$ onto} $\bfm_2$ (for $\bfm_1,\bfm_2 \in \bfM$) \underline{when}
$$
\pi : L_{\bfm_1} \to L_{\bfm_2}
$$
is an order-preserving surjection mapping $\clA_t^{\bfm_1}$ onto 
$\clA_{\pi(t)}^{\bfm_2}$ for each $t \in L_{\bfm_1}$.
\end{definition}

\mn
\begin{remark}
If $t \in L_\bfm$ and $L \subseteq L_\bfm$
, we may abuse notation and write $L_t$ in place of $L \cap L_{\bfm,t}$.

\end{remark}

\mn
\begin{definition}\label{f14}
We define a two-place relation $\leq_\bfM$ (obviously a partial order) on the class of combinatorial templates by:
$$
\bfm_1 \leq_\bfM \bfm_2 \text{ \underline{iff}}
$$
\begin{enumerate}[(a)]
    \item $L_{\bfm_1} \subseteq L_{\bfm_2}$ as linear orders.
\sn
    \item $\bfm_1$ and $\bfm_2$ use the same $\infty$ and the same $t(+)$
    for all $t\in L_{\bfm_1}$.
\sn
    \item $t \in L_{\bfm_1} \Rightarrow \clA_{\bfm_1,t} \subseteq \clA_{\bfm_2,t}$
\end{enumerate}
\end{definition}

\mn
\begin{claim}\label{f17}
$1)$  $\leq_\bfM$ is indeed a partial order on $\bfM$.

\sn
$2)$  If $\LL \bfm_\eps : \eps < \delta\RR$ is $\leq_\bfM$-increasing \underline{then}
$\bfm \defeq \bigcup\limits_{\eps < \delta} \bfm_\eps$ 
(naturally defined) exists, is a $\leq_\bfM$-lub, and is unique.
In particular, 
\begin{enumerate}
    \item [$\boxplus_2$]
    \begin{enumerate}[$(a)$]
        \item $t \in L_\bfm$ \underline{iff} $(\exists \eps < \delta)[t \in L_{\bfm_\eps}]$.
\sn
        \item $s <_\bfm t$ \underline{iff} $(\exists \eps < \delta)[s,t \in L_{\bfm_\eps} \wedge s <_{\bfm_\eps} t]$.
\sn
        \item $\zeta < \delta \wedge t \in L_{\bfm_\zeta} \Rightarrow \clA_{\bfm,t} = \bigcup\limits_{\eps \in [\zeta,\delta)} \clA_{\bfm_\eps,t}$.
    \end{enumerate}
\end{enumerate}

\sn
$3)$ There exists a $\leq_\bfM$-minimal element of $\bfM$. (Specifically, it is an $\bfm$ with $L_\bfm = \varnothing$.)

\sn
$4)$ $(\bfM,\leq_\bfM)$ has amalgamation, and hence the JEP. Moreover, if $\bfm_0 \leq \bfm_\ell$ for $\ell = 1,2$ and (for transparency) $L_{\bfm_0} = L_{\bfm_1} \cap L_{\bfm_2}$, \underline{then} $\bfm_1$ and $\bfm_2$ have a \emph{canonical amalgamation} over $\bfm_0$ --- call it $\bfm_*$.

\noindent
By this we mean 
\begin{enumerate}
    \item [$\boxplus_4$]
    \begin{enumerate}[$(a)$]
        \item $\bfm_* \in \bfM \wedge \bfm_1 \leq_\bfM \bfm_*\wedge \bfm_2 \leq_\bfM\bfm_*$
\sn
        \item $s \in L_{\bfm_*}^+$ \underline{iff} $s \in  L_{\bfm_1}^+ \cup L_{\bfm_2}^+$
\sn
        \item $L_{\bfm_*}^+ \models `r <_\bfm s$\emph{'} \underline{iff}
        \, $`{\bullet_0} \vee {\bullet_1} \vee {\bullet_2}$\emph{,'} where
        \begin{enumerate}
            \item[$\bullet_0$] $L_{\bfm_\ell} \models `s < t$' for some $\ell \in \{0,1\}$.
\sn
            \item  $\big(\exists t \in L_{\bfm_0} \big) \big( \exists k,\ell \in \{1,2\} \big) \big[ r \leq_{\bfm_\ell} t \wedge t \leq_{\bfm_k} s \wedge r \neq s \big]$
\sn
            \item  $r \in L_{\bfm_1}^+ \setminus L_{\bfm_0}^+$, $s \in L_{\bfm_2}^+ \setminus L_{\bfm_0}^+$, and 
            $$
            (\forall t \in L_{\bfm_0}^+)[t <_{\bfm_1} r \wedge t <_{\bfm_2} s].
            $$
        \end{enumerate}
\sn
        \item $\clA_{\bfm,s} = 
        \begin{cases}
            \clA_{\bfm_\ell,s} &\text{if } \ell \in \{1,2\} \wedge s \in L_{\bfm_\ell} \setminus L_{\bfm_0} \\
            \clA_{\bfm_1,s} \cup \clA_{\bfm_2,s} &\text{if } s \in \clA_{\bfm_0}
        \end{cases}$
    \end{enumerate}
\end{enumerate}
\end{claim}

\begin{PROOF}{\ref{f17}}
1) Obvious. 

\sn
2) Straightforward.

\sn
3) Trivial.

\sn
4) We have to check all the clauses in Definition \ref{f11}(1) to verify `$\bfm_* \in \bfM$,' and all clauses in \ref{f17} for $\bfm_\ell \leq_\bfM \bfm_*$.

\bn
\underline{$\bfm_* \in \bfM$:}

\mn
\textbf{Clauses \ref{f11}(1)(a),(b):} Obvious: these were just definitions and notation.



\mn
\textbf{Clause (c):} Also a definition; obvious.



\bn
\underline{$\bfm_\ell \leq_\bfM \bfm_*$} (for $\ell \leq 2$).

Clauses \ref{f14}(a)-(c) are all obvious.
\end{PROOF}

\newpage
\section{Investigating $\bfp \in \bfP$}\label{S2}

\begin{definition}\label{f20}
$\bfP_\bfm$ is the class of $\bfm$-ATIs (see below), and
$$
\bfP \defeq \bigcup\limits_{\bfm \in \bfM} \bfP_\bfm.
$$
\end{definition}
(ATI stands for \emph{abstract template iterations}.)

\mn
\begin{definition}\label{f23}
For $\bfm$ a combinatorial template, we say $\bfp$ is an $\bfm$-\emph{ATI}
\underline{when} it consists of the objects
\begin{enumerate}[$\boxplus$]
    \item
    \begin{enumerate}[(a)]
        \item  $\bfm \in \bfM$ (We may write $L_{\bfp}$ and $L_{\bfp,t}$ instead of $L_{\bfm}$ and $L_{\bfm,t}$, etc.) 
\sn
        \item  Forcing notions $\bbP$ and $\bbP_{\!s}$ for $s \in L_\bfp^+$.
\sn
        \item $\bfo_t \in \bfO_\TR$ (an object, not a name), $\Name\bbQ_t$, and $\name\eta_t$ for $t \in L_\bfp$. We let $S_t = S[t] \defeq S_{\bfo_t}$.
\sn
        \item  $\TR_t = (\TR_{\bfo_t},\leq_{\bfo_t},\plus_{\bfo_t}) = (\TR_t,\leq_t,\plus_t)$ is from $\bfV$.
    \end{enumerate}   
\end{enumerate}       

We demand that the following conditions are satisfied:
\begin{enumerate}
    \item $\prod\limits_t^{<\kappa} \TR_{\bfo_t}$ is a $\kappa$-trunk controller.
\sn
    \item \begin{enumerate}      
        \item $\LL\bbP_{\!s} : s \in L_\bfp^+\RR$ is $\lessdot$-increasing with $s$.
\sn        
        \item Furthermore, we demand $\bbP = \bbP_\infty$.
\sn        
        \item $\name \eta_t$ is a $\bbP_{t(+)}$-name of a member of ${}^{S[t]}2$,\footnote{
            So in the simple case (which we have currently adopted), this will be equal to ${}^{\TR_t}2$.
        } 
        and $\name{\bar\eta} \defeq \LL \name \eta_t : t \in L_\bfp\RR$.
    \end{enumerate}
\sn    
    \item  We require that $p \in \bbP$ \underline{iff}
    \begin{enumerate}
        \item
        $p$ is a function with $\dom(p) \in [L_\bfp]^{< \kappa}$.
\sn
        \item For $s \in \dom(p)$, $p(s)$ is a $\bbP_{\!s}$-name of a member of $\Name\bbQ_s$ of the form $(\tr_{p(s)},\clF_{p(s)})$.
\sn           
        \item  $\tr_{p(s)} \in \TR_t$ is an object, not just a $\bbP_{\!s}$-name.
\sn
        \item $\clF_{p(s)} \defeq \{p_{s,A} : A \in \clA_{p(s)}\}$, where $|\clA_{p(s)}| < \kappa$ 
         But 
        abusing our notation slightly, we allow $A_1 = A_2$ but $p_{s,A_1} \neq p_{s,A_2}$ (but still $|\clF_{p(s)}| < \kappa$  as well). 
\sn
        \item $\gC_{p(s)} \subseteq [\clA_s]^{<\sigma(\tr_{p(s)})}$ (see clause (G)(b) below) and
$$
        p_{s,A} = \bfB(\ldots,\name\eta_{t[j,p_{s,A}]}(\eps^{p_{s,A}}_j),\ldots)_{j < j_{p_{s,A}}}
$$
        \underline{where} $t_j^{p_{s,A}} = t[j,p_{s,A}] \in 
          {}^{ \kappa } A $ 
        and $\bfB = \bfB_{s,A} = \bfB_{p,s,A}$, and
        we make the following demands:
        \begin{enumerate}
            \item $t^{p_{s,A}}_j \in A$,  $\eps^{p_{s,A}}_j \in S_{t_j}$ and  $j_{p_{s,A}} \leq \kappa$.
\sn
            \item $\bfB$ is a $\kappa$-Borel
            function\footnote{
                That is, a function where the pre-image of every element of $S_s$ is a ${\le}\,\kappa$-Borel set. (The point here is absoluteness.)
            }
            from 
            $$
            \prod_{j < j_{p_{s,A}}} S_{t[j,p_{s,A}]}
            $$ 
            to $S_s$ such that the image has cardinality 
            $\leq \kappa$.\\ More concretely, there is 
            (in $\bfV$) an $S'_{p_{s,A}} \in [S_s]^{\le\kappa}$ which contains the image of $\bfB$.
\sn
            \item $\tr_{p_{s,A}} \leq_{\TR_t} \tr_{p(s)}$.
\sn
            \item If  $A \in \clA_{p(s)}$ then $p \rest L_s \Vdash_{\bbP_s} ``\tr_{p(s)} \in 
            {\clE_{p_s,A}}$". 
        \end{enumerate}
\sn
        \item  $S = S_\bfp \defeq \bigcup\limits_{t \in L_\bfp} S_t$.
    \end{enumerate}
\sn
    \item
    \begin{enumerate}
        \item Given $p\in \bbP$ and $s\in \dom(p)$, let $\supp(p(s))$ be the set of all coordinates used in the Borel function $p(s)$ (i.e.~the $t^{p_{s,A}}_j$-s).
        So $\big| \supp(p(s)) \big| \le \kappa$.
\sn
        \item For $p \in \bbP$, we define $\supp(p) \defeq \dom(p) \cup \bigcup\limits_{s\in\dom(p)} \supp(p(s))\in [L_\bfp]^{\leq\kappa}$.
\sn
        \item Note that $\supp(p)\subseteq L_{\bfp,t}$ if 
        $\dom(p)\subseteq L_{\bfp,t}$ (i.e.~if $p\in \bbP_t$).
\sn      
        \item For $L \subseteq L_\bfp$, we set 
        $$
        \bbP_L \defeq \bbP \rest \{p\in\bbP : \supp(p) \subseteq L\}
        $$ 
        and $\bbP_{\!s} \defeq \bbP_{L_{\bfm,s}}.$
\sn      
        \item If $\ell = 1,2$, $s \in I$, $\dom(p^\ell) = \{s\}$, and
        $$
        p^\ell\Big(\bfx_\ell, \big\{\bfB(\ldots,\name\eta_{t[y,(p^\ell)_{s,A_0}]}(\eps[(p^\ell)_{s,A_\ell}]),\ldots)_{i < j_{(p^\ell)_{s,A_\ell}}} \big\} \Big)
        $$
        is as in clause (C)(e) above, \underline{then} there exists some $\kappa$-Borel function 
        $$
        \bfB = \bfB_{\bfx_1,\bfx_2,t}(\ldots,\name\eta_{t_{1,i}}(\eps_{1,i}),\ldots;\ldots , \name\eta_{2,j}(\eps_{2,j}),\ldots )_{\substack{i < \iota_1\\ j < \iota_2}}
        $$ 
        from $\{t_{\ell,i} : \ell \in \{1,2\} \text{ and } i < \iota_\ell\}$ to $\{0,1\}$ which gives the truth value of `$p^1 \leq_{\Name\bbQ_t} p^2$'.
    \end{enumerate}
\sn
    \item For $L \subseteq L_\bfm$ and $p \in \bbP$, we define $p \tser L$ to mean 
        \begin{enumerate}
            \item $\dom(p \tser L) \defeq \dom(p) \cap L$
\sn
            \item For all $s$ in the domain, 
            $$
            (p\tser L)(s) \defeq \big(\tr_{p(s)},\big\{ p_{s,A} : A \in \clA_{p(s)} \cap \clP(L) \big\} \big).
            $$
        \end{enumerate}
\sn
    
\sn
    \item For $t \in L_\bfp$, 
    \begin{enumerate}
        \item $\Name\bbQ_t$ is the {$\bbP_t$-name} of $\Name\bbQ_{\bfo_t}$.
\sn        
        \item If $p \in \bbP_t$ then $\dom(p) \subseteq \supp(p) \in [L_{\bfp,t}]^{\leq\kappa}$.     [This follows.]
\sn            
        \item $\bbP_t$ (and $\bbP_{\!A}$ and $\bbP_L$, for $A \in \clA_\bfp$ and $L \subseteq L_\bfm$) is a $\lepref{\kappa}$-complete $\kappa^+$-cc forcing notion.\footnote{
            The $\kappa^+$-cc follows from {clause (B)(a).}
        } 
\sn            
        \item $\name \eta_t$ is the generic of $\Name\bbQ_t$ (which we may identify with  the subset $\name\eta_t^{-1}(\{1 \}) \subseteq S_t$).
    \end{enumerate}
\sn
    \item 
    \begin{enumerate}
        \item Note that a $\bbP_{\!A}$-generic filter (for 
        $A \in \clA_{\bfm,s}$) lets us evaluate the $\bbP_{r(+)}$-names
        $\name\eta_r$ for $r \in A$, and therefore the value of the 
        $\kappa$-Borel function $p_{s,A}$. 
        
        This way we get a $\bbP_{\!s}$-name for the value, which  we may write as $p_{s,A}[\bfG_{\bbP_{\!A}}]$ or as $p(s)(\name{\bar\eta}\rest A)$.
\sn    
        \item We require that
        $\name \eta_t^{-1}(\{1 \}) = \bigcup\{ \val_t(\tr_{p(t)}) : p \in \Name\bfG_{\bbP_{t(+)}}\}$ (recalling Definition \ref{z20}(1)(F)).
\sn    
        \item If $\kappa \defeq \aleph_0$ then we demand that for every 
        $p \in \bbP_t$ and $m < \kappa$ there exists a $q > p$ such that 
        $|\clF_{q(s)}| \cdot m < \sigma_{\bfo_s}(\tr(q(s)))$ for all $s \in \dom(q)$.
    \end{enumerate}
\sn       
    \item  We require that $p < q$ in\footnote{
        Yes! \underline{Not} `$p \leq q$!'
    } 
    $\bbP$ \underline{iff}
    \begin{enumerate}
        \item $\dom(p) \subseteq \dom(q)$
\sn
        \item If $s \in \dom(p)$ \underline{then} 
        $\tr_{p(s)} \leq_{\bfo_s} \tr_{q(s)}$, $\sigma_{\bfo_s}\!(\tr_{p(s)}) < \sigma_{\bfo_s}\!(\tr_{q(s)})$, and $\clF_{p(s)} \subseteq \clF_{q(s)}$.
    \end{enumerate}
    (Note that for $p \in\bbP$ and $s \in L_\bfp^+$, we have `$p\tser L_{\bfp,s} \in \bbP_{\!s}$' by clause (E).)
        
    Also observe that this is a requirement and \emph{not} a definition, unlike the classical case. 
\sn
    \item When dealing with different ATIs,
        instead of $\bbP$, $\le$, $\bbP_t$, $S_t$, $\Name\bbQ_t$, etc.,
        we may write $\bbP_{\bfp}$, $\le_\bfp$, $\bbP_{\bfp,t}$,
        $S_{\bfp,t}$, $\Name\bbQ_{\bfp,t}$, etc., to indicate that we mean the component of the respective $\bfp$.
\end{enumerate}
\end{definition}

\mn
\begin{remark}\label{f26}
1) Recall that $L_\bfp$ is just a linear order and not necessarily a well-ordering. More concretely,
we do not even exclude the possibility that 
\begin{enumerate}
    \item [$\bullet$] There is an infinite sequence $(s_n)_{n\in\omega}$ with $s_{n+1}\in A_n \in \clA_{s_n}$. 
\end{enumerate}
(unlike \cite{Sh:700}, \cite{Sh:700a}).

\sn 
2) As a consequence: Given $L_\bfm$
and a sequence of (definitions for) $\bbQ_s$ for $s \in L_\bfm$,
it is neither clear that there is an iteration $\bbP$ as above \emph{nor} that it is unique.

\sn 
(In contrast, the usual forcing iteration assumes that 
the index set is well-ordered, and 
we always get a well-defined iteration 
from a sequence of iterands.)

\sn 3)
But if $\bfm$ fails the demand from the bullet in part (1) (so $L_\bfm$ is as in \cite[\S2]{Sh:700}), \underline{then} it is unique.

\sn
4) The following is true:
\begin{enumerate}
    \item [$\circledast$] If $A \in \clA_t$ then $\Vdash_{\bbP_\bfp} ``\name\eta_t$ is generic for $(\bbQ_\bfo,\name\eta_\bfo)$ over $\bfV[\name{\bar\eta} \rest A]$".
\end{enumerate}
\end{remark}

\sn
\begin{claim}\label{f29}
$1)$ If $\bfp \in \bfP$, \underline{then} $\bbP_\bfp$ is $\kappa^+$-cc and  $\lepref{\kappa}$-complete.

\sn
$2)$ In Definition \emph{\ref{f23}}, $\Name\bbQ_t$, $S_t$, and $\TR_t$ are derivable from the other objects.

\sn
$3)$ Any $\leq_{\bbP_\bfp}$-{increasing}\footnote{
    What would occur if we weakened the order by removing the demand `$\sigma_{\bfo_s}\!(\tr_{p(s)}) < \sigma_{\bfo_s}\!(\tr_{q(s)})$' from \ref{f23}(H)(b)?

    \underline{Then} we would have add it as a {caveat} to this statement --- perhaps by stipulating ``$\bar p$ is truly $\bbP_\bfp$-increasing," where
    \emph{truly} means $s \in \dom(p_\eps) \wedge \eps < \zeta < \delta \Rightarrow \sigma(\tr(p_\eps(s))) < \sigma(\tr(p_\zeta(s)))$.
} 
sequence $\bar p = \LL p_\eps : \eps < \delta\RR$ (where $\delta < \kappa$) has a $\leq_{\bbP_\bfp}$-lub $\lim (\bar p)$.
\end{claim}

\begin{PROOF}{\ref{f29}}
1) $\lepref{\kappa}$\underline{-complete:}

Let $\delta < \kappa$ be a limit ordinal and $\bar p = \LL p_i : i < \delta\RR$ be increasing in $\bbP_\bfp$.

Recalling Definition \ref{f23}(H)(b), we know that either $\bar p \rest [j,\delta)$ is constant for some $j < \delta$ \underline{or} $\bar p \rest u$ is $<_{\bbP_\bfp}$-increasing for some unbounded $u \subseteq \delta$. In the first case, $p_j$ is a lub. In the latter case,
\begin{enumerate}
    \item [$\boxplus_1$] For each $\eps \in u$ and $s \in \dom(p_\eps)$, 
    $$
    \Vdash_{\bbP_s}``\LL p_\zeta(s) :\zeta \in u\RR \text{ is {truly} $\leq_{\Name\bbQ_{\bfo_s}}$\!-increasing"}
    $$
    (as defined in the footnote inside clause (3)).
\end{enumerate}
Define the function $p_\delta = \lim(\bar p)$ as follows:
\begin{enumerate}
    \item [$\boxplus_2$]
    \begin{enumerate}
        \item $\dom(p_\delta) \defeq \bigcup\limits_{i<\delta} \dom(p_i)$
\sn
        \item For $s \in \dom(p_\delta)$, let
$$
        p_\delta(s) \defeq \Big( \plus_{\TR_s}( \big\LL \tr_{p_i(s)} : i \in [i_s,\delta) \big\RR), \textstyle\bigcup\limits_{i \in [i_s,\delta)} \clF_{p_i(s)}\Big),
$$
        where $i_s \defeq \min\{i < \delta : s \in \dom(p_i)\}$.
    \end{enumerate}
    
\end{enumerate}
Now check that $p_\delta$ is as required, recalling \ref{f3}(F)(a) and \ref{f23}(H).

\sn
\underline{$\kappa^+$-cc:}

Let $\LL p_\alpha : \alpha < \kappa^+\RR \subseteq \bbP_\bfp$. Without loss of generality, each $p_\alpha$ is like $q$ from \ref{f23}(G)(c) (for $m = 2$). There exists a stationary $S_1 \subseteq S_\kappa^{\kappa^+}$ such that\\ 
$\LL \dom(p_\alpha) : \alpha \in S_1\RR$ is a $\Delta$-system with heart $W_*$. Recalling that $\prod\limits_{s \in W_*} \TR_{\bfo_s}$ 
is a $\kappa$-trunk controller, there exist $\alpha < \beta$ in $S_1$ such that
$$
s \in W_* \Rightarrow \LL \tr_{p_\alpha(s)},\tr_{p_\beta(s)}\RR \in \dom(\plus_{\TR_s}).
$$
If $\kappa > \aleph_0$ then the rest is clear by \ref{f3}(G).

\sn
[Why? Because $\sigma_{\bfo_s}(\tr_{p_\eps(s)})$ is a regular cardinal for some $s \in W_*$.]

\medskip
So assume $\kappa \defeq \aleph_0$. Let $\LL t_\ell : \ell < n\RR$ list $W_*$ in $<_{L_\bfp}$-increasing order, and $t_n \defeq \infty \in L_\bfp^+$. Now fix 
$\alpha < \beta$ from $S_1$, and choose $q_\ell \in \bbP_{t_\ell}$ above 
$$
\{p_\alpha \rest L_{t_\ell},p_\beta \rest L_{t_\ell}\}
$$ 
by induction on $\ell \leq n$. In the induction step, we use our assumption that  
$$
\sigma_{\bfo_s}\!\big( \tr(q_\ell(s)) \big) > 2 \big|\clF_{q_\ell(s)} \big|.
$$

\sn
2) Easy.

\sn
3) Covered in the proof of part (1).
\end{PROOF}

\sn
\begin{definition}\label{f32}
1) We say that $L$ is $\bfp$-closed \underline{when}

\begin{enumerate}
    \item $L \subseteq L_\bfp$
\sn
    \item If $A \in \clA_t$ for some $t \in L$, then $A \subseteq L$.
\end{enumerate}

\sn
2) For 
$L \subseteq L_\bfp$, we define
\begin{enumerate}
    \item $\bbP_L = \bbP_{\bfp,L} \defeq \bbP_\bfp \rest \{p \in \bbP_\bfp : \supp(p) \subseteq L\}$ (as in \ref{f23}(G)(d)).
\sn
    \item $\bfp' = \bfp \rest L$ is defined naturally, with the intention that $\bfp' \in \bfP$. That is, 
    \begin{enumerate}
        \item $\bbP_{\bfp'} \defeq \bbP_{\bfp,L}$
\sn
        \item $\bfo_{\bfp',t} \defeq \bfo_{\bfp,t}$ for all $t \in L$.
\sn
        \item $\clA_{\bfp',t} \defeq \clA_{\bfp,t} \cap \clP(L)$.
    \end{enumerate}
\end{enumerate}

\sn
3) We define ``$\pi$ is an isomorphism from $\bfp'$ onto $\bfp$" naturally. That is,
\begin{enumerate}
    \item $\pi : L_{\bfp'} \to L_\bfp$ is an isomorphism.
\sn
    \item If $t \in L_{\bfp'}$ then $\bfo_{\bfp,\pi(t)} \defeq \bfo_{\bfp',t}$.
\sn
    \item If $t \in L_{\bfp'}$ then $\clA_{\bfp,\pi(t)} \defeq \{\pi[A] : A \in \clA_{\bfp',t}\}$.
\end{enumerate}

\sn
4) We define the partial order $\leq_\bfP$ on $\bfP$ as  $\bfp' \leq_\bfP \bfp$ \underline{iff} $\bfp' = \bfp \rest L_{\bfp'}$. 

\sn
5) If $p \in \bbP_\bfp$ and $L \subseteq L_\bfp$, then $p \tser L$ is defined as follows.
\begin{itemize}
    \item $\dom(p \tser L) \defeq \dom(p) \cap L$
\sn
    \item For all $s$ in the domain, $p\tser L$ sends
    $$
    s \mapsto \big( \tr_{p(s)}, \big\{r \in \clF_{p(s)} : \big(\exists A \in \clA_s \cap \clP(L)\big)\big[ \supp(r) \subseteq A\big] \big\}\big).
    $$
\end{itemize}
\end{definition}

\mn
\begin{observation}\label{f35}\ 
$1)$ If $p\in\bbP$ and $L\subseteq \dom(p)$, then: 
\begin{enumerate}[$(A)$]
    \item $p \tser L$ and $p\rest L$ belong to $\bbP$, with $p \tser L \leq_\bbP p\rest L\le_\bbP p$.
\sn
    \item Moreover, 
    $p\tser L \in \bbP_L$.
\end{enumerate}

\sn   
$2)$  If $L\subseteq L_\bfp$ 
\underline{then}
\begin{enumerate}[$(A)$]
    \item $\bfp \rest L \in \bfP$
\sn
    \item For $p,q$ in $\bbP_L$, we have $p\leq q$ \underline{iff} \emph{\ref{f23}(H)} holds for $\bbP_L$;  i.e.~iff
\begin{itemize}
    \item $\dom(p) \subseteq \dom(q)$ and $\clF_p \subseteq \clF_q$.
\sn            
    \item If $s \in \dom(p)$ and  $r \in \clF_{p(s)}$ with $\supp(r) \subseteq A \in \clA_{\bfp,s}$ (so $A \subseteq L$), \underline{then}
    $$
    q \tser A \Vdash_{\bbP_{\!A}}``\bbQ_{\bfo_t} \models (\tr_p,\{r\}) \leq (\tr_q,\{r\})".
    $$
\end{itemize}
\end{enumerate}


\end{observation}

\begin{PROOF}{\ref{f35}}
Easy.
\end{PROOF}

\mn
\begin{claim}\label{f38}

\sn
$1)$ If $L$ is a 
$\bfp$-closed subset of $L_\bfp$, \underline{then} 
$\bbP_{\bfp,L} \lessdot \bbP_\bfp$.

\sn
$2)$ If $L$ is a 
$\bfp$-closed subset of $L_\bfp$, \underline{then} for all $p \in \bbP_\bfp$, letting $r \defeq p \rest L \in \bbP_{\bfp,L}$, we have
\begin{enumerate}[$\bullet_1$]
    \item $\dom(r) = \dom(p) \cap L$
\sn
    \item If $r \leq q \in \bbP_{\bfp,L}$ \underline{then} 
    $$
    \big( p \rest (\dom(p) \setminus L) \big) \cup q
    $$
    is a lub of $p$ and $q$ (in $\bbP_\bfp$).
\end{enumerate}

\sn
$3)$ If $L_1 \subseteq L_2$ are 
$\bfp$-closed subsets of $L_\bfp$, \underline{then} 
$\bbP_{\bfp,L_2}\, /\, \bbP_{\bfp,L_1}$ is 
$\lepref{\kappa}$-complete and $\kappa^+$-cc.
\end{claim}

\begin{PROOF}{\ref{f38}}
1) Follows by part (2).

\sn
2) Easy.

\sn
3) Straightforward.
\end{PROOF}






\mn
We shall freely use the following claim:
\begin{claim}\label{f40}
Assume $\bfp \in \bfP$ and $\bfp \subseteq \clH(\chi)$, where 
$\chi > \lambda = \lambda^\kappa$. Suppose $N \prec (\clH(\chi),\in)$ is of cardinality $\lambda$ and $[N]^{\leq\kappa} \subseteq N$. \underline{Then} 
$\bfp' \defeq \bfp^N$ (naturally defined) belongs to $\bfP$ and has cardinality $\lambda$.
\end{claim}

\begin{PROOF}{\ref{f40}}
Straightforward.
\end{PROOF}

\mn
\begin{definition}\label{f49}
1) We now define the two-place relation $\leq_\bfP$ as follows.

\sn
Let $\bfp_1 \leq_\bfP \bfp_2$ mean:
\begin{enumerate}[(a)]
    \item $\bfp_\ell$ is an $\bfm_\ell$-ATI for $\ell = 1,2$ (where
    $\bfm_\ell \defeq \bfm_{\bfp_\ell}$; recall that $\bfp_\ell$ determines
    $\bfm_\ell$).
\sn
    \item $\bfm_1 \leq_\bfM \bfm_2$
\sn
    \item For all $t \in L_{\bfm_1}$ we have $\bfo_{\bfp_1,t} = \bfo_{\bfp_2,t}$ (hence $S_{\bfp_2,t} =S_{\bfp_1,t}$ and 
    $\Name\bbQ_t^{\bfp_1} = \Name\bbQ_t^{\bfp_2}$, etc.).
\sn
    \item $\bbP_{\bfp_1} \lessdot \bbP_{\bfp_2}$,
    which implies $\bbP_{\bfp_1,t} \lessdot \bbP_{\bfp_2,t}$ for $t\in L_{\bfm_1}^+$. (This follows by \ref{f38}.)
\sn
    \item \textbf{[Also follows:]} $\Vdash_{\bbP_{\bfp_2}}``\name\eta_t^{\bfp_1} = \name\eta_t^{\bfp_2}$"
    for $t \in L_{\bfm_1}$.
\end{enumerate}

\sn
2)  If $\bfr \leq_\bfP \bfp$ and $p \in \bbP_\bfp$, then we define $p \rest \bfr$ as $p \tser L_\bfr$.

\sn
3) If $\LL \bfp_\alpha : \alpha < \delta\RR$ is $\leq_{\bfP}$-increasing then
``$\bfp \defeq \bigcup\limits_{\alpha< \delta}\bfp_\alpha$" will mean the following:
\begin{enumerate}[(a)]
    \item $\bfp \in \bfP$
\sn
    \item $\bfm_\bfp \defeq \bigcup\limits_{\alpha< \delta} \bfm_{\bfp_\alpha}$
\sn
    \item $\bfp_\alpha \leq_\bfP \bfp$ for all $\alpha < \delta$.
\sn
    \item \textbf{[Follows:]} If $s \in L_{\bfp_\alpha}$ \underline{then}
    $\bfo_{\bfp,s} = \bfo_{\bfp_\alpha,s}$.
\end{enumerate}

\sn
4) We say $\ol \bfp = \LL \bfp_\alpha : \alpha < \alpha_*\RR$ is 
$\leq_{\bfP}$-\emph{increasing continuous} \underline{if} it is 
$\leq_{\bfP}$-increasing and $\bfp_\delta = \bigcup\limits_{\alpha< \delta} \bfp_\alpha$ for every limit $\delta < \alpha_*$.
\end{definition}

\begin{remark}\label{f52}
1) Note that in {parts (3),(4) of} Definition \ref{f49}, for a given $\LL\bfp_\alpha : \alpha < \delta\RR$, it is not \emph{a priori} clear that such $\bfp$ exists --- and even if it does, whether it is unique.

\sn
2) Regarding \ref{f49}(1)(d), does ``$\bbP_{\bfp_1} \lessdot \bbP_{\bfp_2}$"
follow from the rest by \ref{f23}(G)(d) 
by Definition \ref{f14}? {This is not clear.} (See \ref{f26}(2).)

\sn
3) 
We can only show that given $\bfp$ and a $\bfp$-closed
$L \subseteq L_{\bfp}$, we have $(\bfp \rest L) \leq_\bfP \bfp$.
\end{remark}

\sn
\begin{claim}\label{f55}
$1)$ Clause \emph{\ref{f49}(1)(d)} does indeed follow from Claim \emph{\ref{f38}}.

\sn
\emph{1A)} Clause \emph{\ref{f49}(1)(e)} does indeed follow from the rest of \emph{\ref{f49}(1)}.

\sn
$2)$ $\leq_\bfP$ is indeed a partial order on $\bfP$.

\sn
$3)$ If $\LL \bfp_\alpha : \alpha < \delta\RR$ is a $\leq_\bfP$-increasing continuous sequence of $(\partial,\kappa)$-combinatorial templates (Note: when $\kappa > \aleph_0$ this does \underline{NOT} mean that
$\LL \bbP_{\bfp_\alpha} : \alpha < \delta\RR$ is $\subseteq$-increasing continuous!) and $\cf(\delta) \geq \kappa$, \underline{then}
$\bigcup\limits_{\alpha< \delta}\bfp_\alpha$ exists and is unique.

\sn
$4)$ If $\bfp \in \bfP$ and $L$ is $\bfp$-closed, \underline{then} 
$\bfp \rest L \in \bfP$ and $\bfp \rest L \leq_\bfP \bfp$.
\end{claim}

\begin{PROOF}{\ref{f55}}
Straightforward.
\end{PROOF}

\newpage

\section{From $\bfM$ and $\bfP$ to $\bfN$ and $\bfQ$}\label{S3}


For the classes $\bfP$ and $\bfM$, we do not know how to prove the existence of a sufficiently homogeneous member, because when taking unions it is not clear how to define or prove (e.g.) the $\kappa^+$-cc. To solve this, we introduce the classes $\bfQ$ and $\bfN$.

\sn
\begin{definition}\label{d5}
We let $\bfN$ be the class of objects $\bfn$ consisting of\footnote{
    So $\bfm = \bfm_\bfn$, $L_\bfm = L_\bfn$, etc.
}
\begin{enumerate}
    \item $\bfm \in \bfM$
\sn
    \item 
    \begin{enumerate}
        \item $W_\bfn \subseteq L_\bfn$ ($\defeq L_\bfm$)
\sn
        \item If $t \in W_\bfn$ \underline{then}
        \begin{enumerate}
            \item $\lambda_t = \lambda_t^\kappa$
\sn
            \item $\clB_t \defeq \clA_t$ has cardinality $\leq \lambda_{\bfo_t}$. (See \ref{f3}(D)(b).)
\sn
            \item If $A \in \clA_{\bfm,t}$ then $A$ has cardinality 
            $\leq \lambda_{\bfo_t}$ and is a subset of $W_\bfn$.
        \end{enumerate}
    \end{enumerate}
\sn
    \item For each $t \in L_\bfn \setminus W_\bfn$, it also contains 
    $\clB_t = \clB_{\bfn,t}$ and $\clI_t = \clI_{\bfn,t}$ satisfying the following:
    \begin{enumerate}
        \item $\clB_t \subseteq \clA_t$, a family of subsets of $W_{\!\bfn}$.
\sn
        \item If $A \in \clB_t$ then $|A| \leq \lambda_t$.
\sn
        \item $|\clB_{\bfm,t}| \leq 2^{\lambda_{\bfo_t}}$ 
        
        (Really, $\leq \lambda_{\bfo_t}$ would suffice.)
\sn
        \item $\clI_{\bfn,t}$ is a $\kappa$-complete ideal on $\clB_{\bfn,t}$.
\sn
        \item $A \in \clA_{\bfm,t} \Rightarrow \{ B \in \clB_{\bfm,t} : A \cap B \neq \varnothing\} \in \clI_{\bfn,t}$.
    \end{enumerate}
\end{enumerate}
\end{definition}

\mn
\begin{definition}\label{d8}
1) We define a partial order $\leq_\bfN$ on $\bfN$ as follows.

\sn
$\bfn_1 \leq_\bfN \bfn_2$ \underline{iff}
\begin{enumerate}
    \item $\bfm_{\bfn_1} \leq_\bfM \bfm_{\bfn_2}$
\sn
    \item $W_{\bfm_1,s} = W_{\bfm_2,s} \cap L_{\bfm_1}$
\sn
    \item If $s \in L_{\bfm_1}$ then
    \begin{enumerate}
        \item $\lambda_{\bfm_1,s} = \lambda_{\bfm_2,s}$
\sn
        \item $\clB_{\bfm_1,s} = \clB_{\bfm_2,s}$
\sn
        \item $\clA_{\bfm_1,s} \subseteq \clA_{\bfm_2,s}$
\sn
        \item $\clI_{\bfm_1,s} = \clI_{\bfm_2,s} \cap \clP(L_{\bfm_1})$.
    \end{enumerate}
\end{enumerate}

\sn
2) For $\bfn \in \bfN$, we say $L \subseteq L_\bfn$ is $\bfn$-\emph{closed} \underline{when} $A \in \bigcup\limits_{s \in L} \clB_s \Rightarrow A \subseteq L$. 

\sn
3) If $\bfn \in \bfN$ and $L$ is $\bfq$-closed, then $\bfn_\bullet \defeq \bfn \rest L$ is defined naturally: $L_{\bfn_\bullet} \defeq L$, and for all $t \in L$ we have
\begin{itemize}
    \item $\clB_{\bfn_\bullet,t} \defeq \clB_{\bfn,t}$
\sn
    \item $\lambda_{\bfn_\bullet,t} \defeq \lambda_{\bfn,t}$
\sn
    \item $\clA_{\bfn_\bullet,t} \defeq \clA_{\bfn,t} \cap \clP(L)$
\sn
    \item $\clI_{\bfn_\bullet,t} \defeq \clI_{\bfn,t}$.
\end{itemize}
\end{definition}

\mn
\begin{claim}\label{d11}
$1)$ If $\bfn_1 \leq_\bfN \bfn_2$ then $L_{\bfn_1}$ is $\bfn_2$-closed.

\sn
$2)$ If $L$ is $\bfn$-closed then $\bfn \rest L \in \bfN$ and $\bfn \rest L \leq_\bfN \bfn$.

\sn
$3)$ $\leq_\bfN$ is indeed a partial order on $\bfN$.

\sn
$4)$ If $\delta$ is a limit ordinal and $\LL \bfn_\alpha :\alpha < \delta\RR$ is $\leq_\bfN$-increasing, \underline{then} it has a $\leq_\bfN$-least upper bound $\bfn$.

That is:
\begin{enumerate}[$(a)$]
    \item $\bfm_\bfn = \bigcup\limits_{\alpha<\delta} \bfm_{\bfn_\alpha}$ 
\sn
    \item For $\alpha < \delta$, $s \in L_{\bfn_\alpha} \Rightarrow \clB_{\bfn,s} = \clB_{\bfn_\alpha,s}$.
\sn
    \item For $\alpha < \delta$ and $s \in L_{\bfn_\alpha}$, $\clI_{\bfn,s}$ is the closure of $\bigcup\limits_{\beta \in [\alpha,\delta)} \clI_{\bfn_\beta,s}$ under unions of $<\kappa$ members.
\end{enumerate}
\end{claim}

\begin{PROOF}{\ref{d11}}
Straightforward.
\end{PROOF}

\mn
\begin{definition}\label{d20}
We let $\bfQ$ be the class of objects $\bfq$ which consist of
\begin{enumerate}
    \item $\bfn = \bfn_\bfq \in \bfN$
\sn
    \item $\clA_\bfq \defeq \bigcup\limits_{t \in L_\bfn} \clA_{\bfn,t}$
\sn
    \item $\bfp \in \bfP$ such that $\bfm_\bfp = \bfm_\bfn$ and $\lambda_{\bfo_\bfp,t} = \lambda_{\bfn,t}$.
\end{enumerate}
We demand that $\bfq$ satisfies
\begin{enumerate}
    \item [(D)] If $t \in L_\bfq \setminus W_\bfq$ {and} $\bfp_\bullet \in \bfP$ has cardinality $ \leq (|\bbP_\bfp| + \kappa)^\kappa$, {then} the set 
    $$
    \bfY_{\bfp_\bullet,t} \defeq \{A \subseteq W_\bfq \cap L_{\bfq,t} 
    : \bfp_\bfq \rest A \cong \bfp_\bullet\}
    $$
    belongs to $\clI_{\bfq,t}^+$. 

    (If in {(C)} we demand only ``$\clB_{\bfq,t}$ has cardinality $\leq\lambda_{\bfo_t}$," then we should weaken the demand on $A \in \bfY_{\bfp_\bullet,t}$.)
\end{enumerate}
\end{definition}

\mn
\begin{definition}\label{d24}
1) We define a partial order $\leq_\bfQ$ on $\bfQ$ as follows.

\sn
$\bfq_1 \leq_\bfQ \bfq_2$ \underline{iff}
\begin{enumerate}
    \item $\bfn_{\bfq_1} \leq_\bfN \bfn_{\bfq_2}$
\sn
    \item If $s \in L_{\bfq_1}$ \underline{then} $\bfo_{\bfq_1,s} = \bfo_{\bfq_2,s}$.
\sn
    \item $\bbP_{\bfq_1} \lessdot \bbP_{\bfq_2}$. (This follows by \ref{d20}(D).)
\end{enumerate}

\sn
2) For $L$ a $\bfq$-closed (i.e.\ $\bfn_\bfq$-closed) subset of $L_\bfq$, 
let $\bfq \rest L$ be such that $\bbP_{\bfq \rest L} \defeq \bbP_{\bfq,L}$ and $W_{\bfq\rest L} \defeq W_\bfq \cap L$.
\end{definition}

\mn
\begin{claim}\label{d27}
Let $\bfq \in \bfQ$ and $L$ be $\bfq$-closed (that is, $\bfn_\bfq$-closed).

\sn
$1)$ $\bbP_{\bfq,L} = \bbP_{\bfp_\bfq,L} \lessdot \bbP$ and $\bfq \rest L \in \bfQ$.

\sn
$2)$ If $A \in \clA_\bfq$ and $q \in \bbP_\bfq$, \underline{then} there exists $p \in 
\bbP_{\bfq,{A}}$ 
such that \underline{if} 
$\bbP_{\bfq,{A}} 
\models `p \leq r$' \underline{then} $q_* = r \oplus q \in 
\bbP_{\bfq {(+)}}$, 
as defined below, is a common upper bound of $p$ and $q$. 
\begin{itemize}
    \item $\dom(q_*) \defeq \dom(r) \cup \dom(q)$
\sn
    \item $q_*(s) \defeq
    \begin{cases}
        r(s) & \text{if } s \in \dom(r) \setminus \dom(q)\\
        q(s) & \text{if } s \in \dom(q) \setminus \dom(r)\\
        \big( tr(s),\clF_{p(s)} \cup \clF_{q(s)} \big) & \text{if } s \in \dom(q) \cap \dom(r).
    \end{cases}$
\end{itemize}
\end{claim}

\begin{PROOF}{\ref{d27}}
1) Follows from part (2).

\sn
2) By \ref{f40} and Definition \ref{d20}(2).
\end{PROOF}

\newpage

\section{Homogeneity}\label{S4}

\begin{definition}\label{g5}
1) We say $\bfq \in \bfQ_\pre$ \underline{when}
it consists of 
\begin{enumerate}
    \item [(a)] $\bfn_\bfq \in \bfN$\\ (Hence $\bfm_\bfq = \bfm_{\bfn_\bfq} \in \bfM$, $L_\bfq = L_\bfm$, $\clA_\bfq = \clA_{\bfn_\bfq}$, etc. --- see \ref{d20}(A).)
\sn
    \item [(b)] $\overbar\bfo_\bfq = \LL \bfo_{\bfq,t} : t \in L\RR$ as in Definition \ref{f23}.
\end{enumerate}
We demand the following:
\begin{enumerate}
    \item [(c)] For all $A \in \clA_\bfq$ there exists a $\bfq_A \in \bfQ$ with $\bfn_{\bfq_A} = \bfn_\bfq \rest A$ and\\ $\overbar\bfo_{\bfq_A} = \overbar\bfo_\bfq \rest A$. 
\sn
    \item [(d)] If $A \subseteq B$ are members of $\clA_\bfq$ then $\bfq_A = \bfq_B \rest A$.
\end{enumerate}
($\bfq \rest A$ is defined in \ref{f32}(2)(B).)

\sn
2) Let $\leq_\bfQ^\pre$ be the natural order on $\bfQ_\pre$ (see \ref{f32}(2)(B)).

\sn
3) For $\bfq \in \bfQ$, let $\pre(\bfq)$ denote $\bfq$ viewed as a member of $\bfQ_\pre$.
\end{definition}

\mn
\begin{claim}\label{g2}
$1)$ If $\ol\bfq = \LL \bfq_\alpha : \alpha < \delta\RR$ is
$\leq_{\bfQ}$-increasing continuous (see \emph{\ref{f49}(4)})
\underline{then} $\bfq_\delta \defeq \bigcup\limits_{\alpha< \delta}\bfq_\alpha$ exists and is unique, belongs to $\bfQ$, and $\ol\bfq \caret \LL\bfq_\delta\RR$ is $\leq_\bfQ$-increasing continuous.

\sn
$2)$ Every $\bfq \in \bfQ_\pre$ can be expanded to some $\bfq^+ = \bfq(+) \in \bfQ$.  
(By this we mean that $\bfn_{\bfq^+} = \bfn_\bfq$, $\overbar\bfo_{\bfq^+} = \overbar\bfo_\bfq$, $\clA_{\bfq^+} = \clA_\bfq$, and $(\forall A \in \clA_\bfq) \big[ \bfq^+ \rest A = \bfq \rest A \big]$.)


\sn
$3)$ If $\bfr,\bfq \in \bfQ_\pre$ with $\bfr \leq_\bfQ^\pre \bfq$, \underline{then} $\bfr^+ \leq_\bfQ \bfq^+$.
\end{claim}

\begin{PROOF}{\ref{g2}}
1) Let $\bfq_\alpha^* \defeq \pre(\bfq_\alpha)$; that is, a member of $\bfQ_\pre$. 

Clearly,
\begin{enumerate}
    \item [$(*)_1$] $\LL \bfq_\alpha^* : \alpha < \delta\RR$ is $\leq_\bfQ^\pre$-increasing.
\sn
    \item [$(*)_2$] If $I$ is a directed partial order {and} $\overbar\bfq^* = \LL\bfq_s^* : s \in I\RR$ is $\leq_\bfQ^\pre$-increasing continuous, \underline{then} the union $\bfr \defeq \bigcup \overbar\bfq^*$ (naturally defined\footnote{
        In particular:
        \begin{itemize}
            \item $\bfm_\bfr = \bigcup\limits_{s \in I} \bfm_{\bfq_s^*}$ (so $\clA_\bfr = \bigcup\limits_{s \in I} \clA_{\bfq_s^*}$).
\sn
            \item $\bfn_\bfr = \bigcup\limits_{s \in I} \bfn_{\bfq_s^*}$ (so $W_\bfr = \bigcup\limits_{s \in I} W_{\bfq_s^*}$).
        \end{itemize}
    }) 
    is well-defined and is a $\leq_\bfQ^\pre$-lub of the sequence. 
\end{enumerate}
Together with part (2), we are done.

\sn
2) First we will define $\bbP_{\bfq^+}$. Let $\bfn \defeq \bfn_{\bfq^+} = \bfn_\bfq$.
\begin{enumerate}
    \item [$(*)_3$] $p \in \bbP_{\bfq^+}$ \underline{iff}
    \begin{enumerate}
        \item $p$ is a function with domain $\in [L_\bfm]^{<\kappa}$.
\sn
        \item If $s \in \dom(p)$ then $p(s)$ is of the form $(\tr(p(s)),\clF_{p(s)})$, where $\tr(p(s)) \in \TR_{\bfo_s}$ and $|\clF_{p(s)}| < \kappa$.
\sn
        \item If $\name f \in \clF_{p(s)}$, then for some $A \in \clA_{\bfn,s}$ we have $\big(\tr(p(s)),\name f \big) \in \bbP_{\bfq \rest A}$.
    \end{enumerate}
\sn
    \item [$(*)_4$] The order $\leq_{\bbP_{\bfq^+}}$ is defined as follows.

    \sn
    $p \leq_{\bbP_{\bfq(+)}} q$ \underline{iff}
    \begin{enumerate}
        \item $p,q \in \bbP_{\bfq^+}$
\sn
        \item If $A \in \clA_{\bfq^+}$ then $(p \rest A) \leq_{\bbP_{\bfq \rest A}} (q \rest A)$.
    \end{enumerate}
\end{enumerate}
Now check.


\sn
3) Easy.
\end{PROOF}

\mn
\begin{definition}\label{g8}
We define $\paut(\bfq)$ (the \emph{partial automorphisms} of $\bfq$) as the set of $\pi$ such that for some $\bfq$-closed sets $L_1$ and $L_2$, $\pi$ is an isomorphism from $\bfq \rest L_1$ onto $\bfq \rest L_2$.
\end{definition}

\mn
\begin{claim}\label{g11}
$1)$ If $s \in L_\bfq$ then $\bigcup\clB_{\bfq,s} \cup \{s\}$ is a $\bfq$-closed set (recalling \emph{\ref{d8}(2)}).

\sn
$2)$ If $s,t \in L_\bfq$ and $\bfo_s = \bfo_t$, \underline{then} $\{(s,t)\} \in \paut(\bfq)$.

\sn
$3)$ $(\bfQ,\leq_\bfQ)$ has disjoint amalgamation ({as in} \emph{\ref{f17}(4)}): 

\noindent
If $\bfq_1,\bfq_2 \in \bfQ$ and $L_* \defeq L_{\bfq_1} \cap L_{\bfq_2}$ is $\bfq_1$-closed and $\bfq_2$-closed such that $\bfq_1 \rest L_* = \bfq_2 \rest L_*$, \underline{then} there is $\bfp \in \bfQ$ such that $L_\bfp = L_{\bfq_1} \cup L_{\bfq_2}$ and $\bigwedge\limits_{\ell=1,2} \big[\bfq_\ell \leq_\bfQ \bfp \big]$.

\noindent
$4)$ If $\bfo \in \bfO$ \underline{then} there exists $\bfq \in \bfQ$ such that for some $t$, $W_\bfq = \bigcup\clB_{\bfq,t}$, $L_\bfq = W_\bfq \cup \{t\}$, and $\bfo_t = \bfo$.

\sn
$5)$ If $\bfq \in \bfQ$, $\bfo \in \bfO$, $L_*$ is an initial segment of $L_\bfq$, and 
$\clA \subseteq \clP(L_*)$,
\underline{then} there exist $\bfr$ {and $s$} such that:
\begin{enumerate}[$\bullet_1$]
    \item $\bfq \leq_\bfQ \bfr$
\sn
    \item $L_\bfr \setminus L_\bfq = \bigcup\clB_{\bfq,s} \cup \{s\}$
\sn
    \item $\bfo_{\bfr,s} = \bfo$
\sn
    \item If $r \in \bigcup\clB_\bfr \cup \{s\}$ and $t \in L_\bfq$, \underline{then} $L_\bfr \models ``r < t \Leftrightarrow r \in L_*"$.
\sn
    \item $\clA_{\bfr,s} = \clA \cup \clB_{\bfq,s}$.
\end{enumerate}

\sn
$6)$ In part $(5)$, we may allow $\clA \subseteq \clP(L_* \cup \clB_\bfq)$.

\end{claim}

\begin{PROOF}{\ref{g11}}
Straightforward.
\end{PROOF}

\mn
\begin{claim}\label{g14}
Let $\partial = \cf(\partial) \leq \lambda \defeq |\bfO| + \sum\limits_{\TR \in \bfTR} \|\TR\| + \sum\limits_{\bfo \in \bfO} \lambda_\bfo$, where $\lambda_\bfo = \lambda_\bfo^\kappa \geq \|\TR_\bfo\|$ and $\bfO \neq \varnothing$. Suppose $\bfq \in {\bfQ}$ with $\|L_{\bfq}\| \leq \lambda$ such that $L_\bfq \Rightarrow \lambda_{\bfq,s} \leq \lambda_\bfo$.

\underline{Then} for some $\bfr \in \bfQ$:
\begin{enumerate}[$(a)$]
    \item $|L_\bfr| = \lambda$ and $\bfq \leq_\bfQ \bfr$.
\sn
    \item For every $\bfo \in \bfO$, the set $\{s \in L_\bfr : \bfo_{\bfr,s} = \bfo\}$ is dense in $L_\bfr$. 
\sn
    \item There exists an increasing sequence $\olsi L = \LL L_\eps : \eps < \partial\RR$ of subsets of $L_\bfr$, with 
    $\bigcup\olsi L = L_\bfr$, $L_0 \defeq \varnothing$, 
    $L_1 \defeq L_\bfq$, and $\bfr \rest L_1 = \bfq$.
\sn
    \item If $s,t \in L_\bfr$ with $\bfo_{\bfr,s} = \bfo_{\bfr,t} = \bfo$, \underline{then} there is an automorphism $\pi$ of $\bfr$ such that $\pi(s) = t$ and $\pi[L_\eps] = L_\eps$ for every $\eps < \partial$ large enough.
\sn
    \item If $\kappa \leq \theta < \partial$ is such that $\lambda = \lambda^\theta$, $\eps < \partial$, $L_*,L_{**} \in [L_\bfr]^\theta$, and $\pi$ is an isomorphism from $\bfr \rest (L_\eps \cup L_*)$ onto $\bfr \rest (L_\eps \cup L_{**})$ extending $\id_{L_\eps}$,
    \underline{then} there exists an automorphism $\hat\pi : \bfr \to \bfr$ extending $\pi$.
\end{enumerate}
\end{claim}

\begin{PROOF}{\ref{g14}}
We choose $(\bfq_\alpha,\olsi\pi^\alpha)$ by induction on $\alpha < \lambda^+$ such that:
\begin{enumerate}
    \item [$(*)_\alpha$]
    \begin{enumerate}
        \item $\bfq_\alpha \in \bfQ$
\sn
        \item $\bfq_0 \defeq \bfq$ (and $L_{\bfq_0} \defeq \varnothing$ for notational simplicity).
\sn
        \item $\olsi\pi^\alpha = \LL \pi_{\beta,i} : \beta < \alpha,\ i < \lambda \beta +\lambda\RR$
\sn
        \item $\LL \bfq_\beta : \beta \leq \alpha\RR$ is $\leq_\bfQ$-increasing continuous.
\sn
        \item If $\beta < \alpha$ and $i < \lambda\beta + \lambda$, 
        then $\pi_{\beta,i} \in \paut(\bfq_\alpha)$.
\sn
        \item If $\beta < \alpha$ and $i \in [\lambda  \beta,\lambda  \beta +\lambda)$, \underline{then} 
        $\big\LL \pi_{\gamma,i} : \gamma \in [\beta,\alpha) \big\RR$ is $\subseteq$-increasing.
\sn
        \item  If $\gamma + 1 \in [\beta,\alpha)$ and $i \in [\lambda  \beta,\lambda  \beta +\lambda)$, then 
        $$
        L_{\bfq_\gamma} \subseteq \dom(\pi_{\gamma+1,i}) \cap \rang(\pi_{\gamma+1,i}).
        $$
        \item Let $\beta < \alpha$. Each $\pi$ which satisfies 
        `$\bullet_1 \vee \bullet_2$' below will appear somewhere in the sequence $\big\LL \pi_{\beta,i} : i \in [\lambda  \beta,\lambda  \beta +\lambda) \big\RR$.
        \begin{enumerate}
            \item $\pi = \{(s,t)\}$, where $s,t \in L_{\bfq_\beta}$ and $\bfo_{\bfq_\beta,s} = \bfo_{\bfq_\beta,t}$.
\sn
            \item There exist $\gamma < \beta$ and $\bfq_\beta$-closed $L_*$, $L_{**}$ such that $L_* \cup L_{**} \subseteq L_{\bfq_\beta}$, $\lambda = \lambda^{|L_*|}$, and $\pi$ is an isomorphism from $\bfq_\beta \rest (L_{\bfq_\gamma} \cup L_*)$ onto $\bfq_\beta \rest (L_{\bfq_\gamma} \cup L_{**})$ with $\pi \rest L_{\bfq_\gamma} = \id_{L_{\bfq_\gamma}}$.
        \end{enumerate}
    \end{enumerate}
\end{enumerate}

There is no problem in carrying the induction. Now for a club of 
$\delta < \lambda^+$, if $\cf(\delta) = \partial$ then $\bfq_\delta$ is as required.
\end{PROOF}

\bn
\centerline{*\qquad*\qquad*}

\bigskip
We now turn to saccharinity (using a single $\bfo$ for transparency). 

\begin{definition}\label{g19}
Let $\kappa$ be an inaccessible cardinal (or just weakly inaccessible). We say $\bfd$ is a $\kappa$-\emph{parameter} when it consists of
\begin{enumerate}[(a)]
    \item  $\bar\lambda = \LL\lambda_\eta, \lambda_\eta^0, D_\eta : \eta \in \clT\RR$ as in \ref{f9}(4) (so $\clT$ is a subtree of $({}^{\kappa>}\!\kappa,\lhd)$).
\sn
    \item $\bfo_{\bar\lambda}$ as in \ref{f9}(4).
\sn
    \item $\lambda = \cf(\lambda) >\kappa$
\sn
    \item $\bar\rho = \LL \rho_\alpha : \alpha < \lambda\RR \subseteq \prod\limits_{\eta \in \clT} \lambda_\eta$.
\end{enumerate}
Additionally, we demand the following:
\begin{enumerate}
    \item [(e)] \textbf{[Notation:]} For any $L$, let
    \begin{itemize}
        \item $X_{L,\bar\lambda} \defeq $
        \begin{align*}
            \big\{ \bfi = (\zeta,u,g,\nu) : &\ \zeta < \kappa,\ u \subseteq [\kappa]^{<\kappa},\ g : u \to \clT \cap {}^\zeta\kappa,\ \lambda_\nu \ne 1,\\
            & \text{ and }  \eta \in \clT \cap {}^\zeta\kappa \setminus \{\nu\} \Rightarrow \lambda_\eta = 1] \big\}
        \end{align*}
\sn
        \item We define the following partial order on $X_{L,\bar\lambda}$:

        \sn
        $(\zeta_1,u_1,g_1,\nu_1) \leq_{L,\bar\lambda} (\zeta_2,u_2,g_2,\nu_2)$ \underline{iff} 

        $\zeta_1 < \zeta_2$, $\nu_1 \lhd \nu_2$, $u_1 \subseteq u_2$, and $\beta \in u_1 \Rightarrow g_1(\beta) \unlhd g_2(\beta)$.
\sn
        \item Let $Y_{L,\bar\lambda} \defeq \inc_\kappa(X_{L,\bar\lambda})$. 
        
        (That is, $\leq_{L,\bar\lambda}$-increasing sequences from $X_{L,\bar\lambda}$ of length $\kappa$.)
\sn
        \item For $\bar{\bfi} \in Y_{L,\bar\lambda}$, let $E_\bfi \defeq \{\eps < \kappa : \zeta_{\bfi_\eps} = \eps = 
           {\ u_{\bfi_\eps}}\}$;  
        this is a club of $\kappa$.
    \end{itemize}
\sn
    \item [(f)] For any $h : \lambda \to Y_{L,\bar\lambda}$ we can find $(\zeta_*,u_*,g_*,\nu_*) \in X_{L,\bar\lambda}$ and $\mathsf{w} \subseteq \lambda$ such that:
    \begin{itemize}
        \item $\otp(\mathsf{w}) = \lambda_{\nu_*}$
\sn
        \item $(\zeta_{\bfi_\alpha},\nu_{\bfi_\alpha}) = (\zeta_*,\nu_*)$ and 
        $\zeta_* \in E_{h(\alpha)}$ for $\alpha \in \mathsf{w}$.
\sn
        \item $\LL u_{\bfi_\alpha} : \alpha \in \mathsf{w}\RR$ is a $\Delta$-system with heart $u_*$.
\sn
        \item $g_{\bfi_\eps} \rest u_* = g_*$
\sn
        \item $\LL \rho_\alpha(\nu_*) : \alpha \in \mathsf{w}\RR \in {}^\mathsf{w}\lambda_{\nu_*}$ is $<_{L,\bar\lambda}$-increasing. 
    \end{itemize}
\end{enumerate}
\end{definition}

\mn
\begin{claim}\label{g22}
Let $\kappa$ be inaccessible (or weakly inaccessible) such that $\diamondsuit_\kappa$ holds.

\sn
$1)$ If $\kappa^+ = 2^\kappa$ \underline{then} there exists a $\kappa$-parameter $\bfd$.

\sn
$2)$ If $\lambda > \kappa$, $\bbP \defeq \Cohen_\kappa(\lambda)$ 
(the forcing notion adding $\lambda$-many $\kappa$-reals), and $\bar\lambda$, 
$\bfo_{\bar\lambda}$ are as in \emph{\ref{f9}(2)}, then 
$$
\Vdash_\bbP ``\text{there is a $\kappa$-parameter $\bfd$ such that } \bar\lambda_\bfd = \bar\lambda \text{ and } \overbar\bfo_\bfd = \bfo".
$$

\end{claim}

\begin{PROOF}{\ref{g22}}
Easy.
\end{PROOF}

\mn
\begin{claim}\label{g25}
We have  `\,$(A) \Rightarrow (B)$,' where
\begin{enumerate}[$(A)$]
    \item 
    \begin{enumerate}[$(a)$]
        \item $\bfq \in \bfQ$
\sn
        \item $\kappa$ is inaccessible.
\sn
        \item $\bfd$ is a $\kappa$-parameter.
\sn
        \item $\bar\lambda$ is $\kappa$-active (see \emph{\ref{f9}(4)}).
\sn
        \item $s \in L_\bfq \Rightarrow \bfo_{\bfq,s} = \bfo_{\bar\lambda}$
    \end{enumerate}
\sn
    \item If $\bfG \subseteq \bbP_\bfq$ is generic over $\bfV$, \underline{then} in $\bfV[\bfG]$, the set of members $\eta \in \lim(\clT_{\bar\lambda})$ which are generic for $(\bbQ_{\bar\lambda},\name\eta_{\bar\lambda})$ over $\bfV$ is {exactly} $\{\name\eta_{\bfq,s} : s \in L_\bfq\}$.
\end{enumerate}
\end{claim}

\begin{PROOF}{\ref{g25}}
\begin{enumerate}
    \item [$(*)_1$] In $\bfV$ and $\bfV^{\bbP_\bfq}$, we have: `if $\rho \in \prod\limits_{\eta\in\clT_{\bar\lambda}} \lambda_\eta$ \underline{then} $\Name\bfB_\rho \in \id_{\leq\kappa}(\bbQ_{\bar\lambda},\name\eta_{\bar\lambda})$,' where
    $$
    \bfB_\rho \defeq \big\{\nu \in \lim\clT_{\bar\lambda} : (\exists^\kappa \eps < \kappa) \big [\nu(\eps) < \rho(\eps) \big] \big\}.
    $$
\end{enumerate}
This is easily seen to be true. Now assume:
\begin{enumerate}[$(*)_2$]
    \item $\name\eta$ is a $\bbP_\bfq$-name of a member of $\lim \clT_{\bar\lambda}$, and $p \in \bbP_\bfq$ is such that
    $$
    p \Vdash_{\bbP_\bfq} (\forall^\infty \eps < \kappa) \big[\, \name\eta(\eps) > \rho_\alpha(\name\eta \rest \eps) \big]
    $$
    for all $\alpha < \lambda$.
\end{enumerate}

For $\alpha < \lambda$, we can choose $\big\LL(p_{\alpha,\eps},\zeta_{\alpha,\eps},\Lambda_{\alpha,\eps},\nu_{\alpha,\eps}) : \eps < \kappa\RR$ by induction on $\eps$ such that:
\begin{enumerate}
    \item [$(*)_3$]
    \begin{enumerate}
        \item $p_{\alpha,\eps} \in \bbP_\bfq$ is increasing continuously with $\eps$.
\sn
        \item $p \leq p_{\alpha,\eps}$ and $\zeta_* \leq \zeta_{\alpha,\eps}$.
\sn
        \item $(\zeta_{\alpha,\eps},u_{\alpha,
       {\eps}} 
        ,g_{\alpha,
        {\eps}},  
        \nu_{\alpha,\eps}) \in X_{L,\bar\lambda}$, where $\Lambda_{\alpha,\eps} \defeq \rang(g_{\alpha,\zeta})$.

\sn
        \item $\Lambda_{\alpha,\eps} \defeq \big\{ \tr(p_{\alpha,\eps}(\gamma)) : \gamma \in \dom(p_{\alpha,\eps}) \big\}$
\sn
        \item $\LL \zeta_{\alpha,\eps} : \eps < \kappa\RR$ is increasing continuous.
\sn
        \item $p_{\alpha,\eps} \Vdash``\name\eta \rest \zeta_{\alpha,\eps} = \nu_{\alpha,\eps}$"
\sn
        \item $p_{\alpha,0} \Vdash_{\bbP_\bfq}``\text{if $\eps \in [\zeta_{\alpha,0},\kappa)$ then } \name\eta(\eps) > \rho_\alpha (\name\eta \rest \eps)$".
    \end{enumerate}
\end{enumerate}
[Why? For $\eps = 0$ this is easy, by $(*)_2$. For $\eps$ limit use limits, and for $\eps$ a successor ordinal recall that every increasing sequence in $\bbP_\bfq$ of length $<\kappa$ has a lub.]

\begin{enumerate}
    \item [$(*)_4$]
    \begin{enumerate}
        \item Let $h : \lambda \to Y_{L_\bfq,\bar\lambda}$ be defined by 
        $\beta \mapsto \LL \bar\bfi_\alpha : \alpha < \beta\RR$.
\sn
        \item Applying \ref{g19}(f) to our choice of $h$, we can choose $\mathsf{w}$ and $(\zeta_*,u_*,g_*,\nu_*)$ as there.
    \end{enumerate}
\end{enumerate}

Lastly, we define $q$.
\begin{enumerate}
    \item $\dom(q) \defeq \bigcup\limits_{\alpha \in \mathsf{w}}\dom(p_{\alpha,\zeta_*})$
\sn
    \item $q(s) \defeq
    \begin{cases}
        p_\alpha(s) &\text{if } s \in \dom(p_{\alpha,\zeta_*}) \setminus u_*\\
        \big( \tr(p_{\alpha_*}(s)), \bigcup\limits_{\alpha \in \mathsf{w}} \clF_{p_\alpha(s)} \big) &\text{if } s \in u_*,
    \end{cases}$

    \noindent
    where we fix $\alpha_*$ as (e.g.) $\min u_*$.
\end{enumerate}
[Why is $q \in \bbP_\bfq$? Because if $\alpha \in \mathsf{w}$ and $\varrho \in \clT \cap \bigcup\limits_{\eps \le \zeta_*}{}^\eps\!\kappa$, then $\varrho \neq \nu_* \Rightarrow \lambda_{\nu_*} < \lambda_\varrho^0$.]

Also, $q$ is a $\leq_{\bbP_\bfq}$-upper bound of $\{p_{\alpha,\zeta_*} : \alpha \in \mathsf{w}\}$, but this implies 
$$
q \Vdash `\alpha \in \mathsf{w} \Rightarrow \name \eta(\zeta_*) > \rho_\alpha(\zeta_*)\text{'.}
$$ 
As $\LL \rho_\alpha(\zeta_*) : \alpha < \lambda_{\nu_*}\RR \in {}^{\lambda_{\nu_*}}\lambda_{\nu_*}$

is increasing, we get a contradiction to 
$p_\alpha \Vdash `\name\eta(\zeta_*) < \lambda_{\name\eta \rest \zeta_*}$' and\\ $p_\alpha \Vdash ` \lambda_{\name\eta \rest \zeta_*} = \lambda_{\nu_*}$'.
\end{PROOF}

\mn
\begin{fact}\label{g20}
The conclusion of \ref{g25} holds for any inaccessible $\kappa$ (allowing a preliminary forcing) without any extra set-theoretic assumptions.
\end{fact}

\begin{PROOF}{\ref{g20}}
If we force by $\bbP \defeq \mathsf{Levy}(\kappa^+,2^\kappa)$ \underline{then} 
$$
\Vdash_\bbP``2^\kappa = \kappa^+ \text{ and the desired conclusion holds"},
$$ 
so we have $\LL\name p_\alpha : \alpha \in \Name S\RR$: names of the desired objects. But $\bbP_\bfq = \bbP_\bfq^{\bfV[\bbP]}$, so we can easily finish.
\end{PROOF}

\mn
As in \ref{g25} and \ref{g20}, 
\begin{claim}\label{g23}
We have `$(A) \wedge (B) \Rightarrow (C)$', where
\begin{enumerate}[$(A)$]
    \item  $\kappa \defeq \aleph_0 < \lambda = \lambda^{\aleph_0} \leq 2^\kappa$ 
\sn
    \item 
    \begin{enumerate}[$(a)$]
        \item $\bfq \in \bfQ$
\sn
        \item  $\bfo = \bfo_\bfn^2$ is as in \emph{\ref{f9}(1)}.
\sn
        \item $s \in L_\bfq \Rightarrow \bfo_{\bfq,s} = \bfo$.
    \end{enumerate}
\sn
    \item  If $\bfG \subseteq \bbP_\bfq$ is generic over $\bfV$, \underline{then} in $\bfV[\bfG]$, the set of members $\eta \in \lim(\clT_{\bar\lambda})$ which are generic for $(\bbQ_\bfn^2,\name\eta_\bfn^2)$ over $\bfV$ is {exactly} $\{\name\eta_{\bfq,s} : s \in L_\bfq\}$.
\end{enumerate}
\end{claim}

\begin{PROOF}{\ref{g23}}
Essentially by \cite[\S6, Claim 21]{Sh:1067}. There the iteration is FS, but this does not cause any serious problem.
Also, we can prove this similarly to the claims above.
\end{PROOF}

\newpage

\section{Variants, and the theorem}\label{S5}

\begin{discussion}\label{i2}
1) In \S1-4 we chose the simplest version of the iteration $\bfq$: specifically, $\TR_\bfq$ was simple. Here we will present a more general version.

\sn
2) 
Considering a $\leq_{\bbP_\bfq}$-increasing sequence $\LL p_i : i < \delta\RR$, where $\delta < \kappa$ is a limit ordinal, we do not always have a least upper bound. This occurs when $\cf(\delta) \in \Theta$: otherwise we would have to add the existence of lubs in the demands in \ref{f3}.

\sn
3) Also, \ref{f29}(3) has to be changed.

\sn
4)  We have the option of generalizing `$\kappa$-trunk controller' to `$(\kappa,\Theta,\Upsilon)$-trunk controller,' but to save on parameters we will not pursue this.
(But we will remark when the proof requires changes.)

\sn
5) In the iteration, we did not deal with the case where $\Name\bbQ_{\bfp,t}$ is defined from {a} parameter which is a $\bbP_{\bfp,t}$-name (even in this case it is a $\bbP_{\bfp,t}$-name of a real).
\end{discussion}

\mn
\begin{fact}\label{i5}
We can repeat {\S1-4} with the following changes:

\mn
$1)$ In {\ref{f2}(3),(4)}, 
\begin{enumerate}[$\bullet_1$]
    \item Every $\TR \in \bfTR$ is a $(\kappa,\Theta,\Upsilon)$-trunk controller (see {\ref{z20}(2)}) so $\bfTR$ is not necessarily simple.
\sn
    \item Again, we demand that the $<\kappa$-support product of a sequence of members of $\bfTR$ is always a $\kappa$-trunk controller. (This is not crucial.)
\sn
    \item $\theta = \cf(\theta) > \kappa$ and $\alpha < \kappa \Rightarrow |\alpha|^\kappa < \theta$
\end{enumerate}

\mn
$2)$ In {\ref{f3}(A)}, 
$$
\TR_\bfo = (|\TR_\bfo|,\leq_\bfo,\plus_\bfo,S_\bfo,\val_\bfo).
$$

\mn
$3)$ In \ref{f3}(F)(b):
\begin{itemize}
    \item In the conclusion, generally only an upper bound exists.
\sn
    \item But if $\cf(\delta) \in \Theta$ then the full conclusion holds.
\end{itemize}

\mn
$4)$ In \ref{f7}(2), we demand only
\begin{itemize}
    \item Any truly increasing sequence in $\bbQ_\bfo$ of length $<\kappa$ has an upper bound.
\sn
    \item Any truly increasing sequence of length $\delta <\kappa$ with $\cf(\delta) \in \Theta$ has a least upper bound.
\end{itemize}

\mn
$5)$ We replace {\ref{f11}(1)(c)} by
\begin{enumerate}
    \item [$(c)'$] $\clA_t \subseteq [L_{\bfq,t}]^{<\theta}$
\end{enumerate}

\mn
$6)$ In Definition \ref{f14}, we omit clause (f).

\mn
$7)$ In {\ref{f29}(3)}, if $\cf(\delta) \in \Theta$ {then} yes, {there are} lubs. (Otherwise, we just have upper bounds.)

\mn
8) In $(*)_{s,i}$ in the proof of \ref{f38}(2), we replace clause (c)
with
\begin{enumerate}
    \item [(c)$'$] $A_{s,i}^* \cap A_{s,j}^* = \varnothing$ for all $j < i$, and $A_{s,i}^* \cap A = \varnothing$ for all $A \in \clA_{p(s),i}$.
\end{enumerate}

\mn
9) Nevertheless, if we weaken clause (1)\,$\bullet_2$ \underline{then} {in} \ref{g14}(b) we have to add the assumption that $\TR_\bfo \times \prod\limits_{t \in L_\bfq}^{<\kappa} \TR_{\bfo_t} \notin \bfO^+$ is a $\kappa$-trunk controller.
\end{fact}

\mn
\begin{theorem}\label{i8}
Adopting the context of \emph{\ref{i2}}, we have $(A) \Rightarrow (B)$ and\\ $(A) \wedge (B) \Rightarrow (C)$, where
\begin{enumerate}[$(A)$]
    \item
    \begin{enumerate}[$(a)$]
        \item We have $`\bullet_1 \vee \bullet_2$,\emph{'} where
        \begin{enumerate}
            \item $\kappa$ is inaccessible and $\bfo \defeq \bfo_{\bar\lambda}$ (so $\bbQ_\bfo = \bbQ_{\bar\lambda}$).
\sn
            \item $\kappa \defeq \aleph_0$ and $\bfo \defeq \bfo_\kappa^2$ as in \emph{\ref{f9}(1)}.
        \end{enumerate}
\sn
        \item $\lambda > \theta = \cf(\theta) > \kappa$
    \end{enumerate}
\sn
    \item There exist $\overbar\bfq$ and $\bar t$ such that
    \begin{enumerate}[$(a)$]
        \item $\overbar\bfq = \LL \bfq_\alpha : \alpha \leq \theta\RR \subseteq \bfQ$. Let $\bfq \defeq \bfq_\theta$.
\sn
        \item $|L_{\bfq_\alpha}| = \lambda$
\sn
        \item $\bar t = \LL t_\alpha : \alpha < \theta\RR$
\sn
        \item $t_\alpha \in L_{\bfq_{\alpha+1}}$ and $L_{\bfq_\alpha} \in \clA_{\bfq,t_\alpha}$
\sn
        \item $\bfo_{\bfq,t} \defeq \bfo$
\sn
        \item For all $r,s \in L_\bfq$ there exists a $\pi \in \aut(L_\bfq)$ such that $\pi(r) = s$.
\sn
        \item Moreover, if $s <_{L_\bfq} t_1$ and $s <_{L_\bfq} t_2$, \underline{then} there exists a $\pi \in \aut(L_\bfq)$ such that $\pi(t_1) = t_2$ \emph{and} $\pi \rest L_{\leq s}$ is the identity.
    \end{enumerate}
\sn
    \item Letting $\Name X \defeq \{\name\eta_{\bfq,s} : s \in L_\bfq\}$ be a $\bbP_\bfq$-name, we have that 
    $$
    \bfV_{\!1} \defeq \HOD(\bfV \cup \Name X)^{\bfV^{\bbP_\bfq}}
    $$ 
    satisfies the following:
    \begin{enumerate}[$(a)$]
        \item It is a model of $\ZF$.
\sn
        \item $\id_{<\theta}(\bbQ,\name\eta)$ has measurability.
\sn
        \item $\bfV_{\!1}$ has $\DC_{<\theta}$.
    \end{enumerate}
\end{enumerate}
\end{theorem}

\begin{PROOF}{\ref{i8}}
Straightforward.
\end{PROOF}

\mn
\begin{remark}\label{i11}
We may consider an alternative to \ref{i5}, retaining ``$\clA_{\bfm,t}$ is a collection of pairwise disjoint sets." \underline{But} when defining $\bfm \leq_\bfM \bfn$, we replace ``$\clA_{\bfm,s} \subseteq \clA_{\bfn,s}$" by
\begin{itemize}
    \item For every $A \in \clA_{\bfm,s}^+$ there exists a \underline{uni}q\underline{ue} $B \in \clA_{\bfn,s}^+$ such that $B \supseteq A$.
\end{itemize}
\end{remark}

\bibliographystyle{amsalpha}
\bibliography{shlhetal}
\end{document}